\documentclass[10pt,a4paper,reqno]{amsart} 

\usepackage{amsmath}
\usepackage{bbm}
\usepackage{mathtools}
\usepackage{amssymb}
\usepackage{graphicx}
\usepackage{color}
\usepackage{latexsym}  
\usepackage{amssymb} 
\usepackage[mathcal]{eucal}
\usepackage{fancybox}
\usepackage{algorithm,algorithmicx,algpseudocode}
\usepackage{xcolor}
\usepackage{rotating}
\usepackage[a4paper]{hyperref}

\algrenewcommand\algorithmicindent{2.5em}
\algrenewcommand{\algorithmiccomment}[1]{{\color{brown} \hfill$\blacktriangleright$ #1 }}

\DeclareMathOperator{\Act}{Act}

\DeclareMathOperator{\cc}{cc}
\DeclareMathOperator{\cycl}{cycl}

\newcommand{\fig}[3]{
\begin{figure}[h!]
\begin{center}
 \includegraphics #1
 \end{center}
\vspace{-5pt}
\caption{ #2}
\label{#3}
\end{figure}
}

\newcommand{\N}{\mathbb N}

\newcommand{\diffs}[2]{#1 \, \ominus \, #2}

\newcommand{\tut}{\mathcal T_{G,\Delta}}
\newcommand{\sing}[1]{\{#1\}}
\newcommand{\ens}[1]{\left\{#1\right\}}
\newcommand{\ar}[1]{\{#1,#1'\}}
\newcommand{\contract}[2]{\cccc(#1,#2)}
\newcommand{\delete}[2]{\dddd(#1,#2)}

\definecolor{darkgray}{gray}{0.4}

\DeclareMathOperator{\cccc}{contract}
\DeclareMathOperator{\dddd}{delete}

\newcommand{\bSi}{\textbf{S}$_i$}
\newcommand{\mSi}{\textbf{S}_i}
\newcommand{\bSe}{\textbf{S}$_e$}
\newcommand{\mSe}{\textbf{S}_e}
\newcommand{\bI}{\textbf{I}}
\newcommand{\bL}{\textbf{L}}
\newcommand{\ua}[1]{\underset{#1}{\rightarrow}}

\newtheorem{prop}{Proposition}[section]

\newtheorem{lem}[prop]{Lemma}
\newtheorem{defi}[prop]{Definition}

\newtheorem{theo}[prop]{Theorem}

\newtheorem{conjecture}{Conjecture}
\newtheorem{cor}[prop]{Corollary}

\catcode`\@=11
\def\section{\@startsection{section}{1}%
 \z@{.7\linespacing\@plus\linespacing}{.5\linespacing}%
 {\normalfont\bfseries\scshape\centering}}

\def\subsection{\@startsection{subsection}{2}%
  \z@{.5\linespacing\@plus\linespacing}{.5\linespacing}%
  {\normalfont\bfseries\scshape}}

\def\subsubsection{\@startsection{subsubsection}{3}%
 \z@{.5\linespacing\@plus\linespacing}{-.5em}
  {\normalfont\bfseries\itshape}}
\catcode`\@=12

\title{A general notion of activity for the Tutte polynomial} 
\author{Julien COURTIEL}

\begin{document}

\maketitle

\begin{abstract}
In the literature can be found several descriptions of the Tutte polynomial of graphs. Tutte defined it thanks to a notion of \textit{activity} based on an ordering of the edges. Thereafter, Bernardi gave a non-equivalent notion of the activity where the graph is embedded in a surface. In this paper, we see that other notions of activity can thus be imagined and they can all be embodied in a same notion, the \textit{$\Delta$-activity}. We develop a short theory which sheds light on the connections between the different expressions of the Tutte polynomial.
\end{abstract}

\section{Introduction}

Intended as a generalization of the chromatic polynomial \cite{whitney,tutte54}, the Tutte polynomial is a graph invariant playing a fundamental role in graph theory.
 This polynomial 
is essentially the generating function of the spanning subgraphs, counted by their number of connected components and cycles. Because of its wide study and its universality property of deletion/contraction \cite{broxley}, the Tutte polynomial is often used as reference polynomial when it comes to interlink graph polynomials from different research fields. For example, the Potts model in statistical physics \cite{fk}, the weight enumerator polynomial in coding theory \cite{curtis}, the reliability polynomial in network theory \cite{oxley-welsh} and the Jones polynomial of an alternating node in knot theory \cite{jonespoly} can all be expressed as (partial) evaluations of the Tutte polynomial. The reader can be referred for more information to the survey \cite{chap-tutte} written by Ellis-Monaghan and Merino.

It is well-known that the Tutte polynomial can be also defined as the generating function of the spanning trees $T$ counted according to the numbers $i(T)$ and $e(T)$ of internal and external active edges:
\begin{equation}\label{first}
T_G(x,y) = \sum_{T\textrm{ spanning tree}} x^{i(T)}y^{e(T)}.
\end{equation}
Such a description appeared for the first time in the founding paper of Tutte \cite{tutte54}. Tutte's notion of activity required to linearly order the edge set. More recently, Bernardi \cite{bernardi-tutte} gave a new notion of activity which was this time based on an embedding of the graph. 
The two notions are not equivalent, but they both satisfy \eqref{first}. 
One can also cite the notion of external activity introduced by Gessel and Sagan \cite{GesselSagan} involving the depth-first search algorithm and requiring a linear ordering on the vertex set.

The purpose of this paper is to unify all these notions of activity.
We thereby define a new notion of activity, called \textit{$\Delta$-activity}. Its definition is  based on a new combinatorial object named \textit{decision tree}, to which the letter $\Delta$ refers.  We show that each of the previous activities is a particular case of $\Delta$-activity. Moreover, we see that the $\Delta$-activity enjoys  most of the properties that were true for the other activities, like Crapo's property of partition of the subgraphs into intervals \cite{crapo}.

Here is an overview of the paper. We begin in Section \ref{DNM} by summarizing the definitions and the notations thereafter needed. For instance an \textit{activity}  denotes a function that maps a spanning tree onto a set of active edges. An activity is said to be \textit{Tutte-descriptive} if it describes the Tutte polynomial in the sense of \eqref{first}. Section \ref{sec:activity} outlines four families of Tutte-descriptive activities, including the aforementioned three ones. A new fourth family of activities is described, named \textit{blossoming activities}. It is based on the transformation of a planar map into blossoming trees \cite{Sch97,bdg2002}. In the context of a study on maps, it constitutes an interesting alternative to Bernardi's notion of activity.
 
 Section \ref{s:alg} introduces the notion of $\Delta$-activity through an algorithm.  This first definition of the $\Delta$-activity is strongly reminiscent of what we could see in \cite{gordon-traldi}. In this paper, Gordon and Traldi stated that the Tutte's active edges for a subset $S$ 
of a matroid $M$ correspond to the elements that are deleted or contracted as an isthmus or a loop during the so-called \textit{resolution}\footnote{A resolution of a matroid is a sequence of deletions and contractions which reduces the matroid into the empty matroid.} of $M$ with respect to $S$. 
Our approach is very similar, except that we do not consider our ``resolution" in a fixed order as Gordon and Traldi did. We also prove that the $\Delta$-activities are Tutte-descriptive.
In Section \ref{s:ord}, we state some properties connecting $\Delta$-activity and edge ordering. For instance, we see that an edge is active when it is maximal (for a certain order depending on the spanning tree) inside its fundamental cycle/cocycle. 
In Section \ref{s:partition} we set forth a partition of the set of subgraphs that is the counterpart of Crapo's partition \cite{crapo} for \mbox{$\Delta$-activities}. This partition results from the equivalence relation naturally induced by the algorithm of Section \ref{s:alg}. We deduce from this an enlightening proof of the equivalence between the principal descriptions of the Tutte polynomial. 

In Section \ref{sec:spec}, we show that the four families of activities of Section \ref{sec:activity} can be defined in terms of $\Delta$-activities. We thus prove in an alternative way that these activities are all Tutte-descriptive. Some extra properties can be deduced using the theory from Section \ref{s:alg}, \ref{s:ord}, \ref{s:partition}. In Section \ref{s:com}, we discuss about some extensions of the $\Delta$-activities, as the (easy) generalization to the matroids. We conclude the paper by a conjecture that would emphasize the relevance of the $\Delta$-activities: In rough terms it states that the ``interesting" activities exactly coincide with the $\Delta$-activities.

\section{Definitions, notations, motivations}
\label{DNM}

\subsection{Sets and graphs}

\subsubsection{Sets}

The set of non-negative integers is denoted by $\N$. We denote by $|A|$ the cardinality of any set $A$. When we say that a set $A$ is the disjoint union of some subsets $A_1,\dots,A_k$, this means that $A$ is the union of these subsets and that $A_1,\dots,A_k$ are pairwise disjoint. We then write $A = \biguplus_{i=1}^k A_i$. For any pair of sets $A$, $B$, we denote by $\diffs A B$ the symmetric difference of $A$ and $B$.\footnote{We do not use the notation $A \, \triangle \, B$ because the triangle $\triangle$ can be easily mistaken for the letter $\Delta$, very used in this paper.} Let us recall that the symmetric difference is commutative and associative. 

\subsubsection{Graphs, subgraphs, intervals} \label{sss:interval}
In this paper, the graphs will be finite and undirected. Moreover, they may contain loops and multiple edges. For a graph $G$, the set of vertices is denoted by $V(G)$ and the set of edges by $E(G)$.

A \textit{spanning subgraph} of $G$ is a graph $S$ such that $V(S)=V(G)$ and $E(S) \subseteq E(G)$. Unless otherwise indicated, all subgraphs will be spanning in this paper. A subgraph is completely determined by its edge set, therefore we identify the subgraph with its edge set. For instance, given $A$ a set of edges, we will allow ourselves to write $S \subseteq A$ if $S$ is a subgraph of $G$ only made of edges that belong to $A$. For any subgraph $S$, we denote by $\overline S$ the complement subgraph of $S$ in $G$, i.e. the subgraph $\overline S$ such that $V(\overline S) = V(G)$ and $E(\overline S) = E(G) \backslash E(S)$.
If $S$ is a subgraph of $G$ and $A \subseteq E(S)$, we write $S \backslash A$ to denote the subgraph of $G$ with edge set $E(S) \backslash A$. 
Given a subgraph $S$ of $G$, an edge $e \in E(G)$ is said to be \emph{internal} if $e \in S$, \emph{external} otherwise.

The collection of all subgraphs of $G$ can be ordered via inclusion to obtain a boolean lattice. A \textit{subgraph interval} $I$ denotes an interval for this lattice, meaning that there exist two subgraphs $H^-$ and $H^+$ such that $I$ is the set of subgraphs $S$ satisfying $H^- \subseteq S \subseteq H^+$. In this case, we write $I = [H^-,H^+]$.

\subsubsection{Cycles and cocycles}

A \emph{path} is an alternating sequence $(v_1,e_1,v_2,\dots,e_k,v_{k+1})$ of vertices and edges such that the endpoints of $e_j$ are $v_j$ and $v_{j+1}$ for all $j \in \ens{1,\dots,k}$. A path with no repeated vertices and edges, with the possible exception $v_1=v_{k+1}$, is called a \emph{simple} path. A path is said to be \emph{closed} when  $v_1 = v_{k+1}$. A \emph{cycle} is the set of edges of a simple and closed path. For instance, a loop can be seen as a cycle reduced to a singleton. A graph (or a subgraph) with no cycle is said to be \emph{acyclic}.

A \emph{cut} of a graph $G$ is a set of edges $K$ such that the endpoints of each edge of $K$ are in two distinct connected components of $\overline K$. A \emph{cocycle} is a cut which is minimal for inclusion. (Therefore the deletion of a cocycle exactly increases the number of connected components by one.) An \emph{isthmus} is an edge whose deletion increases the number of connected components. (In other terms, an isthmus is a cocycle reduced to a singleton.) Note that a subgraph $S$ of $G$ is connected if and only if the set $\overline S$ does not contain any cocycle of $G$.

An edge which is neither an isthmus, nor a loop, is said to be \emph{standard}.

\noindent \textbf{Example. }Consider the graph of Figure \ref{exg}. The set $\ens{a,c,d}$ is a cycle but $\ens{a,c,d,g,i,j}$ is not (because the associated path is not simple). The set $\ens{b,h}$ is a cocycle but  $\ens{b,h,i,j}$ is not.

\fig{[scale=1.2]{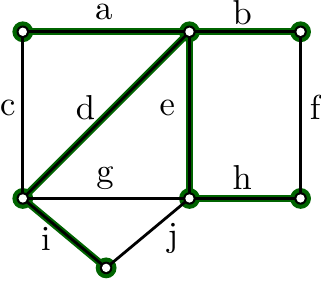}}{A graph equipped with a spanning tree in bold edges}{exg}

A \emph{spanning tree} is a (spanning) subgraph that is a tree, that is to say a connected and acyclic graph. Fix $T$ a spanning tree of a graph $G$. The \emph{fundamental cycle} of an external edge $e$ is the only cycle contained in the subgraph $T \cup \{e\}$, i.e. the cycle made of $e$ and the unique path in $T$ linking the endpoints of $e$. Similarly, the \emph{fundamental cocycle} of an internal edge $e$ is the only cocycle contained in $\overline T \cup \sing e$, i.e. the cocycle made of the edges having exactly one endpoint in each of the two subtrees obtained from $T$ by removing $e$.

\noindent \textbf{Example. }Consider the spanning tree from  Figure \ref{exg}. The fundamental cycle of $j$ is $\ens{d,e,i,j}$, the fundamental cocycle of $e$ is $\ens{e,f,g,j}$.

\subsubsection{Embeddings}

We now define combinatorial map as Robert Cori and Toni Machi did in \cite{cori-these,cori-machi}. A \emph{map} $M = (H,\sigma,\alpha)$ is a finite set of half-edges $H$, a permutation $\sigma$ of $H$ and an involution without fixed point $\alpha$ on $H$ such that the group generated by $\sigma$ and $\alpha$ acts transitively on $H$. A map is \textit{rooted} when one of its half-edges, called the \textit{root}, is distinguished. In this paper, \emph{all the maps are rooted}.

For every map $M = (H,\sigma,\alpha)$, we define its \emph{underlying graph} as follows: We form a vertex for each cycle of $\sigma$, and an edge for each cycle of $\alpha$. A vertex is incident with an edge if the corresponding cycles have a non-empty intersection. Observe that such a graph is always connected since $\sigma$ and $\alpha$ act transitively on $H$.

\begin{figure}[h!]
\begin{center}
\includegraphics[scale=1]{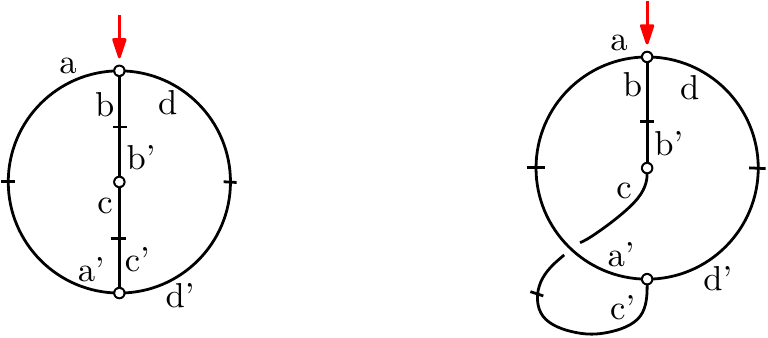}
\end{center}
\caption{Two different embeddings of the same graph.}
\label{fig:embe}
\end{figure}

A \emph{combinatorial embedding} of a connected graph $G$ is a map $M$ such that the underlying graph of $M$ is isomorphic to $G$. Every embedding of $G$ is in correspondence with a \emph{rotation system} of $G$, that is to say a choice of a circular ordering of the half-edges around each vertex of $G$. The rotation system actually corresponds to the above permutation $\sigma$. The permutation $\alpha$ is automatically given by the matching of the half-edges. When an embedding of $G$ is considered, we write the edges of $G$ as pairs of half-edges (for instance $e = \ar h$). 

A map $M = (H,\sigma,\alpha)$ can be graphically represented in the following way: We draw the underlying graph of $M$ in such a manner that the cyclic ordering of the half-edges in counterclockwise order around each vertex matches the corresponding cycle in $\sigma$. If we want to indicate the root, we add an arrow just before\footnote{in counterclockwise order} the root half-edge, pointing to the root vertex. In other terms, if $h_0$ denotes the root, we put an arrow between the half-edges $h_0$ and $\sigma^{-1}(h_0)$.

For example, Figure \ref{fig:embe} shows two different maps with the same underlying graph. The left map corresponds to $(H,\sigma,\alpha)$ and the right map to $(H,\sigma',\alpha)$, where $H = \ens{a,a',b,b',c,c',d,d'}$,  $\sigma = (a\,b\,d)(a'\,d'\,c')(b'\,c)$, $\sigma' = (a\,b\,d)(a'\,c'\,d')(b'\,c)$ and \mbox{$\alpha = (a\,a') (b\,b') (c\,c') (d\,d')$} (the permutations are written in cyclic notation). Each of these maps is rooted on the half-edge $a$.

Given a map $M = (H,\sigma,\alpha)$, we say that a half-edge $h_2$ \emph{immediately follows} another half-edge $h_1$ when $h_2 = \sigma \circ \alpha (h_1)$. In Figure \ref{fig:embe}, the half-edge that immediately follows $a$ for the left map is $d'$, while it is $c'$ for the right map.

\subsubsection{Edge contraction and deletion}

\emph{Edge deletion} is the operation that removes an edge $e$ from the edge set of a graph $G$ but leaves its endpoints unchanged. The resulting graph is denoted by $\delete G e$.

\emph{Edge contraction} is the operation that removes an edge from a graph by merging the two vertices it previously connected. More precisely, given an edge $e$ with endpoints $v$ and $w$ in a graph $G$, the contraction of $e$ yields a new graph $G'$ where $w$ has been removed from the vertex set and any edge which was incident to $w$ in $G$ is now incident to $v$ in $G'$.
The resulting graph $G'$ is denoted by $\contract G e$. 

We extend the operations of deletion and contraction to maps, as follows. Let us consider $M = (H,\sigma,\alpha)$ an embedding of a graph $G$ and $e = \ens{h_1,h_2}$ an edge of $G$. Let $H'$ be $H \backslash \ens{h_1,h_2}$ and $\alpha'$ the involution $\alpha$ restricted to $H'$. If $e$ is not an isthmus, we define $\delete M e$ as the map $(H',\sigma_d,\alpha')$ where
\begin{equation}
\sigma_d(h) = \left\{ 
\begin{array}{ll} 
\sigma \circ \sigma \circ \sigma(h) & \textrm{if }(\sigma(h) = h_1\textrm{ and }\sigma(h_1)=h_2)\textrm{ or } (\sigma(h) = h_2\textrm{ and }\sigma(h_2)=h_1), \\
\sigma \circ \sigma(h) & \textrm{if }(\sigma(h) = h_1\textrm{ and }\sigma(h_1) \neq h_2)\textrm{ or } (\sigma(h) = h_2\textrm{ and }\sigma(h_2) \neq h_1), \\
\sigma(h) & \textrm{otherwise}. \\
\end{array}
\right. 
\end{equation}
Similarly, if $e$ is not a loop, we define $\contract M e$ as the map $(H',\sigma_c,\alpha')$ where
\begin{equation}
\sigma_c(h) = \left\{ 
\begin{array}{ll} 
\sigma \circ \sigma (h) & \textrm{if }(\sigma(h) = h_1\textrm{ and }\sigma(h_2)=h_2)\textrm{ or } (\sigma(h) = h_2\textrm{ and }\sigma(h_1)=h_1), \\
\sigma \circ \alpha \circ \sigma (h) & \textrm{if }(\sigma(h) = h_1\textrm{ and }\sigma(h_2) \neq h_2)\textrm{ or } (\sigma(h) = h_2\textrm{ and }\sigma(h_1) \neq h_1), \\
\sigma(h) & \textrm{otherwise}. \\
\end{array}
\right. 
\end{equation}
We do not want to define the deletion of an isthmus of the contraction of a loop since it could bring about the disconnection of the map\footnote{Even for a contraction. Take for instance the map $(H,\sigma,\alpha)$ with $H = \ens{a,a',b,b',c,c'}$, $\sigma = (a\,b\,a'\,c) (b')(c')$ and $\alpha = (a\,a') (b\,b') (c\, c')$. If $e = \ar a$, we can check that $\sigma_c = (b)(b')(c)(c')$.}.

One can check that $\delete M e$ (resp. $\contract M e$) is indeed an embedding of $\delete G e$ (resp. $\contract G e$). If $h_1$ is the root of the map $M$, we choose $\sigma_d(h_1)$ (resp. $\sigma_c(h_1)$) as the root of $\delete M e$ (resp. $\contract M e$).


%
%
%
%
%
%

\subsection{The Tutte polynomial}

\subsubsection{Definition}
Two parameters are important to define the Tutte polynomial. The first one is the \emph{number of connected components} of a subgraph $S$, denoted by $\cc(S)$. Recall that by convention each subgraph is spanning. This implies that the subgraph of $G$ with no edge has $|V(G)|$ connected components. The second parameter is the \emph{cyclomatic number} of $S$, denoted by $\cycl(S)$. It can be defined as 
\begin{equation} \cycl(S) = \cc(S) + |S| - |V(G)|. \end{equation}
It equals the minimal number of edges that we need to remove from $S$ to obtain an acyclic graph. In particular, $\cycl(S)=0$ if and only if $S$ is a forest.

\begin{defi} The \emph{Tutte polynomial} of a graph $G$ is 
\begin{equation} T_G(x,y) = \sum_{S\textrm{ subgraph of }G}(x-1)^{\cc(S)-\cc(G)}(y-1)^{\cycl(S)},
\label{eq:deftutte}
\end{equation}
where $\cc(S)$ and $\cycl(S)$ respectively denote the number of connected components of $S$ and the cyclomatic number of $S$.
\end{defi}

For example, consider the graph of Figure \ref{fig:extut}. Let us list all its subgraphs with their contributions to the Tutte polynomial: there are one subgraph with no edge (contribution \mbox{$(x-1)^2$}), four subgraphs with one edge (contribution $4\,(x-1)$), five acyclic subgraphs with two edges (contribution $5$), the subgraph $\ens{a,d}$ (contribution $(x-1)\,(y-1)$), four subgraphs with three edges (contribution $4\,(y-1)$) and the whole graph (contribution $(y-1)^2$). Thus the Tutte polynomial of this graph is $$(x-1)^2 + 4 \, (x-1) +  5 + (x-1)\,(y-1) + 4\,(y-1) + (y-1)^2,$$ which can be rewritten as $x^2 + x + x\,y + y + y^2$.

\subsubsection{Properties} One can easily duduce from \eqref{eq:deftutte} that if the graph $G$ is the disjoint union of two graphs $G = G_1 \uplus G_2$, then the Tutte polynomial of $G$ is the product of the two other Tutte polynomials: \mbox{$T_G(x,y) = T_{G_1}(x,y) \times T_{G_2}(x,y)$}. We can \textit{de facto} restrict our study to connected graphs. \emph{From now, all the graphs $G$ we consider are connected.}
We also assume that our graphs have at least one edge. (Otherwise, the Tutte polynomial is equal to $1$.)  

Let us  recall the relations of induction satisfied by the Tutte polynomial, due to Tutte himself \cite{tutte54}.

\begin{prop} Let $G$ be a graph and $e$ be one of its edges. The Tutte polynomial of $G$ satisfies:
\begin{equation}
T_G(x,y) = \left\{ \begin{array}{ll} T_{\contract G e}(x,y) + T_{\delete G e}(x,y) & \textrm{if }e\textrm{ is standard,} \\
x\,T_{\contract G e}(x,y) & \textrm{if }e\textrm{ is an isthmus,} \\
y\, T_{\delete G e}(x,y) & \textrm{if }e\textrm{ is a loop.} 
 \end{array} \right.
\label{eq:ind}
\end{equation}
\end{prop} 

Since the Tutte polynomial of a graph with one edge is equal to $x$ or $y$, the previous proposition implies by induction the following property:

\begin{cor}The Tutte polynomial of any graph has non-negative integer coefficients in $x$ and $y$.
\end{cor} 

It is natural to ask if a combinatorial interpretation exists for these coefficients. Tutte found in 1954 a characterization of his polynomial in terms of an "activity" based on a total ordering of the edges. Some decades later, Bernardi gave a similar characterization with a notion of activity this time related to an embedding of the graph. The precise definitions will be given in Section~\ref{sec:activity}.

\subsubsection{Activities}

Let us formalize the notion of activity.  (The following definitions are not conventional, as the activity of a spanning tree usually denotes the number of active edges.)

An \emph{activity} is a function that maps spanning trees of $G$ on subsets of $E(G)$.
We say that an activity $\psi$ is \emph{Tutte-descriptive} if the Tutte polynomial of $G$ is equal to 
\begin{equation}
T_G(x,y) = \sum_{T\textrm{ spanning tree of }G} x^{|\mathcal I(T)|}\,y^{|\mathcal E(T)|}, 
\label{eq:compact} 
\end{equation}
where $\mathcal I (T)= \psi(T) \cap T$ and $\mathcal E(T) = \psi(T) \cap \overline T$.

An \emph{internal activity} (resp. \emph{external activity}) is a function that maps any spanning tree $T$ of $G$ onto a subset of $T$ (resp. a subset of $\overline T$). We say that an internal activity $\mathcal I$ (resp. external activity $\mathcal E$) can be \emph{extended into a Tutte-descriptive activity} if there exists a Tutte-descriptive activity $\psi$ such that $\psi(T) \cap T = \mathcal I(T)$ (resp. $\psi(T) \cap \overline T = \mathcal E(T)$) for any spanning tree $T$ of $G$. 

The objective of this paper is to introduce several families of Tutte-descriptive activities, and to describe a general framework from which every of these activities can be deduced.  


\section{Four families of activities}
\label{sec:activity}

In this section, we introduce four families of activities, one of which is new. All these families are (pairwise) not equivalent, meaning that no family is included in another. We will prove in Section \ref{sec:spec} that any of these activities is Tutte-descriptive, or can be extended into Tutte-descriptive activities.

\subsection{Ordering activity (Tutte)} 
\label{ss:ord}

The first family of activity was defined by Tutte in \cite{tutte54}, as follows.

Consider $G$ a graph. We equip it with a linear ordering on the edge set. An external (resp. internal) edge of a spanning tree $T$ is said to be \emph{ordering-active} if it is minimal in its fundamental cycle (resp. cocycle). The \emph{ordering activity} is the function that sends every spanning tree $T$ onto the set of its ordering-active edges. This activity naturally depends on the chosen ordering of the edges.

\begin{figure}[h!]
\begin{center}
\includegraphics[scale=1]{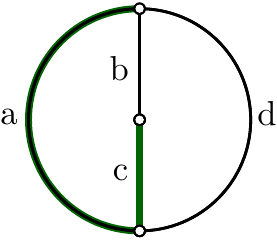}
\end{center}
\caption{A graph with spanning tree $T = \sing {a,c}$.}
\label{fig:extut}
\end{figure}

\noindent \textbf{Example:} Consider the graph of Figure \ref{fig:extut}. The edges are ordered alphabetically, that is to say $a < b < c < d$. With this ordering, the spanning tree  $T = \sing {a,c}$ induces only one internal active edge, $a$, and no external active edge. Indeed, $b$ is not externally active since its fundamental cycle is $\{a,b,c\}$. Similarly, $c$ is not internally active because $b$ belongs to its fundamental cocycle. 

Tutte established that the ordering activity is always Tutte-descriptive. In particular, the sum \eqref{eq:compact}, where $\mathcal I(T)$ and $\mathcal E(T)$ denote the sets of internal and external ordering-active edges of a spanning tree $T$, does not depend on the chosen linear ordering, although the activity clearly does.

\subsection{Embedding activity (Bernardi)}
\label{ss:bernardi}

Bernardi defined in \cite{bernardi-tutte} other activities that are well adapted to the notion of maps.

We consider $M_G = (H,\sigma,\alpha)$ an embedding of a graph $G$, that we root on a half-edge denoted by $h_0$. To each spanning tree $T$, we associate a \emph{motion function} $t$ on the set $H$ of half-edges by setting
\begin{equation} 
\label{motionfunction}
t(h) = \left\{ \begin{array}{ll} \sigma(h) & \textrm{if }h\textrm{ is external,} \\ \sigma \circ \alpha (h)  & \textrm{if }h\textrm{ is internal.}    \end{array} \right. 
\end{equation}
If the notation is ambiguous, we will write $t(h;T)$.

\begin{figure}[h!]
\begin{center}
\includegraphics[scale=1]{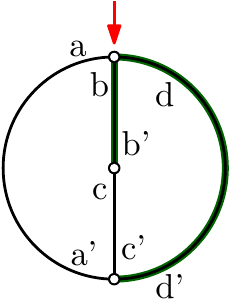} \hspace{3cm} \includegraphics[scale=1]{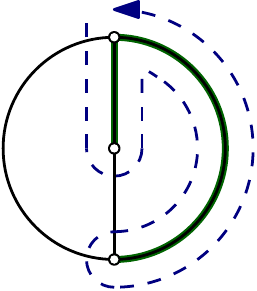}
\end{center}

\caption{Example of an embedded graph equipped with a spanning tree and representation of the corresponding tour.}
\label{fig:exber}
\end{figure}

The motion function characterizes the \emph{tour} of a spanning tree. In informal terms, a tour is a counterclockwise walk around the tree that follows internal edges and crosses external edges (see Figure \ref{fig:exber}). Bernardi proved the following result \cite{bernardi-tutte}.
\begin{prop} 
\label{tourcyclique}
For each spanning tree, the motion function is a cyclic permutation of the half-edges.
\end{prop}
Thus, we can define a linear order on the set $H$ of half-edges, called the ($M_G$,$T$)-order, by setting 
$h_0 < t(h_0) < t^2(h_0) < \dots < t^{|H| - 1}(h_0)$, where $h_0$ is the root of $M_G$. This order can be then transposed on the set of edges: we say that $e = \ar{h_1} < e' = \ar{h_2}$ when $\min(h_1,h_1') < \min(h_2,h'_2)$.

Then, an external (resp. internal) edge is said to be $(M_G,T)$-\emph{active} if it is minimal for the $(M_G,T)$-order in its fundamental cycle (resp. cocycle). The \emph{embedding activity} is the function that associates with a spanning tree $T$ the set of $(M_G,T)$-active edges of $G$.

\noindent \textbf{Example:} Take  the embedded graph from Figure \ref{fig:exber}, that we will denote $M_G$, rooted on $a$ and equipped with the spanning tree $T=\sing {\ar b,\ar d}$. The motion function for this spanning tree is the cycle $(a,b,c,b',d,c',a',d')$.  So the half-edges are sorted for the $(M_G,T)$-order as follows: $a < b < c < b' < d < c' < a' < d'$. Thus, the $(M_G,T)$-order for the edges is $\ar a  < \ar b  < \ar c < \ar d $. There is one external active edge, $\ar a $, and one internal active edge, $\ar b$. The edge $\ar c$ (resp. $\ar d$) is not active since $\ar b$ (resp. $\ar a$) is in its fundamental cycle (resp. cocycle). 

It was proven by Bernardi that any embedding activity is Tutte-descriptive, whatever the chosen embedding is.

\subsection{Blossoming activity}
\label{ss:blo}
We are going to give a new family of activities, named \emph{blossoming activities}. We first define it for internal edges only.
As for any embedding activity, we need beforehand to embed the graph $G$ and root it. We denote by $M$ the resulting  map. 

Given a spanning forest $F$, Algorithm \ref{prune} outputs a spanning tree denoted $\tau(F)$. \\

\begin{algorithm}[h!]
\caption{Computing $\tau(F)$}
\label{prune}
\begin{algorithmic}
\Require $F$ spanning forest of $M$.
\Ensure $\tau(F)$ spanning tree of $M$.
\State $h \leftarrow$ root of $M$;
\State $M' \leftarrow M$;
\While {some edges are not visited} 
	\State $e \leftarrow$ edge that contains $h$;
	\State $h \leftarrow$ half-edge that immediately follows $h$ in $M'$; \Comment{i.e. $h \leftarrow \sigma \circ \alpha(h)$}
	\If{$e$ is not an isthmus of $M'$ \textbf{and} $e \notin F$}
		\State $M' \leftarrow \delete {M'} {e};$
	\EndIf 
\EndWhile
\State $\tau(F) \leftarrow M'$
\State \Return $\tau(F)$;
\end{algorithmic}
\end{algorithm}

\noindent \textbf{Informal description:} Starting from the root, we turn counterclockwise around the map. \emph{After} we walk along an edge, we remove it whenever the deletion of this edge leaves the map connected (i.e. it is not an isthmus) \emph{and} it is external. We stop the algorithm when every edge has been visited. \\

\begin{figure}[h!]
\begin{center}
\includegraphics[scale=0.9]{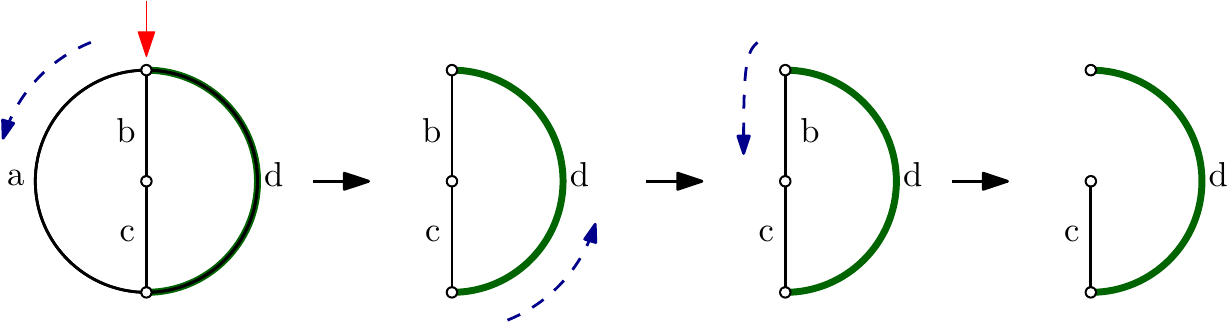}
\end{center}
\caption{Run of Algorithm \ref{prune} with spanning forest $F = \sing {d}$.}
\label{fig:exblo}
\end{figure}

\noindent \textbf{Example:} We consider the first embedded graph of Figure \ref{fig:exblo} with spanning forest $F = \sing {d}$. The run of Algorithm \ref{prune} with this input is illustrated on the same figure. We have $\tau(F) = \ens{c,d}$.

Let us justify the termination of the algorithm.

\begin{prop}
For any spanning forest $F$, Algorithm \ref{prune} with input $F$ stops and outputs a spanning tree.
\end{prop}

\begin{proof} Consider a step in the run of Algorithm \ref{prune} where the map $M'$ is not a tree. Then some edges incident to the face\footnote{The faces of a combinatorial map $(H,\sigma,\alpha)$ are the cycles of $\sigma \circ \alpha$.} that contains the half-edge $h$ in $M'$ form a cycle. This cycle contains an external edge because $F$ is acyclic. So the face that contains $h$ includes an external non-isthmus edge. As $e$ takes for successive values in Algorithm \ref{prune} the edges of the face incident to $h$, the edge $e$ will be eventually both external and a non-isthmus edge. At this time, $e$ will be removed from $M'$. 

In other terms, while $M'$ is not a tree, an edge is deleted. Since $M'$ remains connected, $M'$ will end as a tree. Then, $M'$ has only one face and every edge has been visited by the algorithm.\end{proof}

\noindent \textbf{Remark 1.} For those who are familiar with \cite{Sch97,bdg2002}, Algorithm \ref{prune} 
is related to the transformation of a map into a \emph{blossoming tree}. If we do not take into account the so-called \textit{buds} and \textit{leaves}, the blossoming tree that corresponds to the map $M$ is exactly $\tau(\emptyset)$. Subtler connections exist and are described in the author's PhD \cite{courtiel-these}.



Given a spanning tree $T$ of $G$, we say that an internal edge $e$ is \emph{blossoming-active} if
$$\tau(T \backslash e) = T.$$
The \emph{internal blossoming activity} is the function that maps a spanning tree $T$ onto the set of its internal blossoming-active edges. 

\noindent \textbf{Example.} For the first map of Figure \ref{fig:exblo} with spanning tree $T = \sing {c,d}$, the only internally blossoming-active edge is $c$. Indeed, we saw that $\tau(\ens{d})=T$ and  we observe that $\tau(\ens{c})= \ens{b,c}$.

We will see with Proposition \ref{prop:blo} that the internal blossoming activity can be extended into a Tutte-descriptive activity, for any embedding of the graph. Observe that the blossoming-active edges are not defined in terms of their fundamental cocycles, but such a description exists, as is shown in Corollary \ref{corblo}.

\noindent \textbf{Remark 2.} Let us give an alternative definition of this activity, based on an assignment of charges, when $M$ is planar \footnote{No proof will be given, it would be a bit off topic. However the reader can refer to \cite[p. 214]{courtiel-these}.}. A \textit{charge} is an element from $\ens{-1,+1}$. Given a spanning tree $T$ of $M$, we run Algorithm \ref{prune} and whenever we delete an edge $e$, we add a charge $-1$ to the departure vertex and a charge $+1$ to the arrival vertex. At the end of the algorithm, there only remains a rooted plane tree, namely $T$, where each vertex has several charges. Let $e$ be an internal edge. Deleting $e$ splits $T$ into two components. The one which does \textit{not} contain the root is the \textit{subtree} corresponding to $e$.
 
 \begin{prop} \label{p:newcar}
   Given any spanning tree, if an internal edge is blossoming-active, then the sum of charges in the corresponding subtree is $0$ or $1$. 
If $M$ is  planar, then the converse is true.  
  \end{prop}

\begin{figure}[h!]
\begin{center}
\includegraphics[scale=0.9]{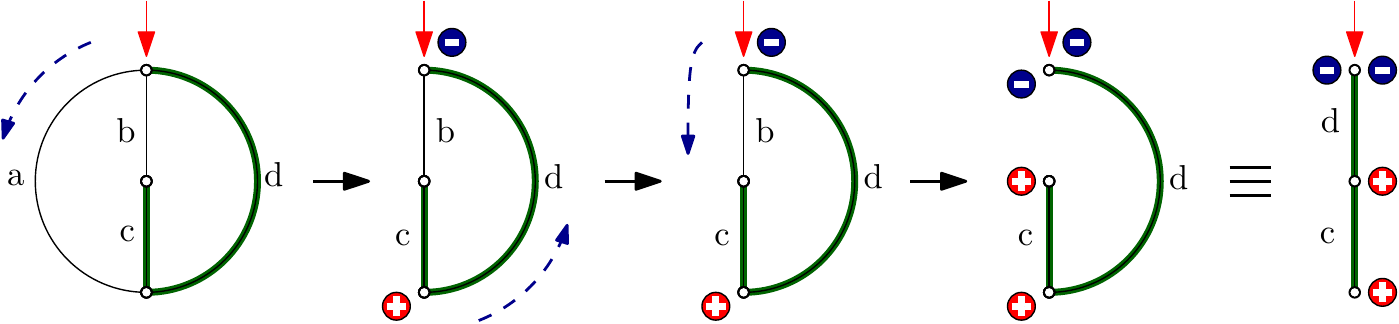}
\end{center}
\caption{Run of the version with charges of Algorithm \ref{prune} with spanning tree $T = \ens {c,d}$.}
\label{fig:excharge}
\end{figure}

An example is shown in Figure \ref{fig:excharge}. Four charges have been distributed on the plane tree $T$.
If we delete $d$, the corresponding subtree has two vertices with one charge $+1$ each, so the charge of the subtree is $+2$. Hence the edge $d$ is not blossoming-active. On the contrary, if we delete $c$, the subtree has only one charge $+1$: the edge $c$ is blossoming-active. We have checked on this example Proposition \ref{p:newcar}: the only internal blossoming-active edge is $c$. 

\begin{figure}[h!]
\begin{center}
\includegraphics[scale=0.9]{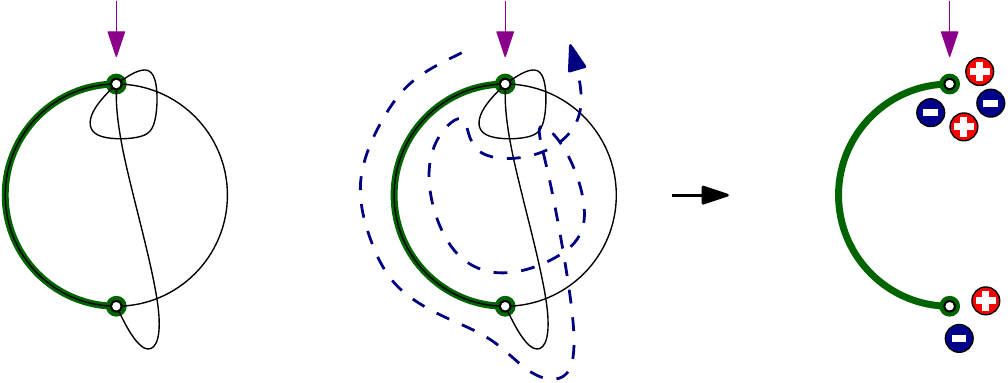}
\end{center}
\caption{A counterexample for Proposition \ref{p:newcar} with $M$ not planar.}
\label{cexba}
\end{figure}

When $M$ is not planar, the previous property is false. A counterexample is shown in Figure \ref{cexba}. The subtree corresponding to the only internal edge has charge $0$ but this edge is not blossoming active. (It will be deleted at the first step of Algorithm \ref{prune} if we remove it from the spanning tree.)

\subsection{DFS activity (Gessel, Sagan)}
\label{ss:dfs}

Gessel and Sagan described in \cite{GesselSagan} a notion of external edge activity for external edges based on the Depth First Search (DFS in abbreviate). 

Consider a graph $G$. We assume that $V(G) = \ens{1,\dots,n}.$ Moreover, we assume that $G$ does not have multiple edges. This allows us to avoid technical details, which can be easily included if needed. Every edge will be denoted by the pair of integers that corresponds to its endpoints, for example $e=\ens{1,2}$.

Given a (not necessarily connected) graph $H$, Algorithm \ref{DFS} computes the \textit{(greatest-neighbor) DFS forest} of $H$, denoted by $\mathcal F(H)$.

\begin{algorithm}[h!]
\caption{DFS forest of a graph}
\label{DFS}
\begin{algorithmic}[5]
\Require graph $H$.
\Ensure $\mathcal F(H)$ spanning forest of $H$.
\State $\mathcal F(H) \leftarrow \emptyset$;
\While {there is a unvisited vertex} 
	\State \textbf{mark} the least unvisited vertex of $H$ \textbf{as visited}; \label{linedfs}
	\While {there is a visited vertex with unvisited neighbors} 
		\State $v \leftarrow$ the most recently visited such vertex;
		\While{$v$ has a unvisited neighbor} 
			\State $u \leftarrow$ the greatest unvisited neighbor of $v$;
			\State $e \leftarrow \ens{u,v}$;
			\State \textbf{mark} $u$ \textbf{as visited};
			\State $v \leftarrow u$; 
			\State \textbf{add} $e$ in $\mathcal F(H)$;
		\EndWhile
	\EndWhile 
\EndWhile
\State \Return $\mathcal F(H)$
\end{algorithmic}
\end{algorithm}

\noindent \textbf{Informal description:}  We begin by the least vertex. We proceed to the DFS of the graph $H$ that favors the largest neighbors. Each time we move from a vertex to another, we add the edge that we have followed to the DFS forest. When we have visited all the vertices of a connected component, we reset the process, starting from the least unvisited vertex. 

\begin{figure}[h!]
\begin{center}
\includegraphics[scale=1]{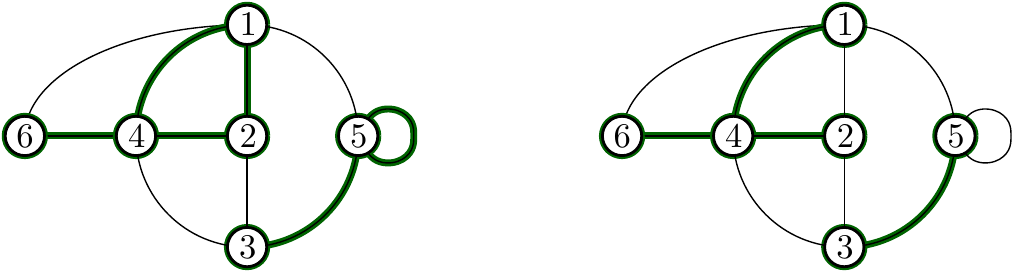}
\end{center}
\caption{Left: a subgraph indicated by thick green edges. Right: its DFS forest.}
\label{fig:dfs}
\end{figure}

Here we only apply the algorithm to a subgraph $H$ of $G$. An example is shown in Figure \ref{fig:dfs}: during the DFS of the subgraph on the left, the vertices will be visited in the order $1,4,6,2,3,5$. The resulting DFS forest is represented on the right. Note that $S$ and $\mathcal F (S)$ have the same number of connected components. 

Given a spanning forest $F$ of $G$, we say that an external edge is \emph{DFS-active} if $$\mathcal F\left(F \cup \ens{e}\right) = F.$$
The \textit{external DFS activity} is the function that sends every spanning forest onto the set of its external DFS-active edges.

For instance, let us go back the spanning forest $F$ of Figure \ref{fig:dfs} (right). It has two external DFS-active edges: $\ens{1,2}$ and $\ens{5,5}$. On the contrary, the edge $e=\ens{1,6}$ is not DFS-active since $\mathcal F\left(F \cup \ens{e}\right)$ is equal to \mbox{$\ens{\ens{1,6},\ens{6,4},\ens{4,2},\ens{3,5}}$}.

There is a natural question about Algorithm \ref{DFS}: given a spanning forest $F$ of $G$, can we describe the set of subgraphs $S$ such that $\mathcal F(S)$ equals $F$? The notion of external DFS activity answers to this question. Indeed, for any subgraph $S$ and for any spanning forest $F$ of $G$, we have $\mathcal F(S) = F$ if and only if \mbox{$S \in [F, F \cup \mathcal E(F)]$}, where $\mathcal E(F)$ denotes the set of external DFS-active edges (see Subsection \ref{sss:interval} for the definition of intervals).

Let us cite this alternative characterization of the externally DFS-active edges (Lemma 3.2 from \cite{GesselSagan}).
\begin{prop} \label{lemGS}
Let $F$ be a spanning forest of $G$. An external edge $e$ is DFS-active if and only if:
\begin{itemize}
\item $e$ is a loop, or
\item $e = \ens{u,v}$ where $v$ is a descendant\footnote{that is to say a vertex in the same component of $u$, visited after $u$} of $u$, and $w > v$, where $w$ is the child of $u$ on the unique path in $F$ linking $u$ and $v$. (We also say that $(w,v)$ is an inversion.)
\end{itemize}
The last case is depicted in Figure \ref{descendant}.
\end{prop}

\fig{[scale = 1.5]{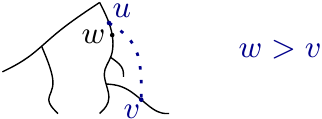}}{A non-loop DFS-active external edge.}{descendant}

Gessel and Sagan proved that the Tutte polynomial of $G$ satisfies:
\begin{equation} \label{DFStut}
T_G(x,y) = \sum_{F\textrm{ spanning forest of }G} (x-1)^{\cc(F)-1} y^{|\mathcal E(F)|},
\end{equation}
where $\mathcal E$ denotes the external DFS activity. We are going to show that the external DFS activity, restricted to spanning trees, can be extended into a Tutte-descriptive activity (cf. Prop. \ref{dfstut}). 

In the same paper, Gessel and Sagan defined a notion of \emph{external activity with respect to NFS} (NFS for Neighbors-first search). We will not detail it but this activity also falls within the scope of this work.
Also note that a notion of edge activity based on Breadth-First Search is conceivable.

\section{Algorithmic definition of $\Delta$-activity}
\label{s:alg}
We are going to introduce a meta-family of Tutte-descriptive activities, named $\Delta$-activities. This family includes all the Tutte-descriptive activities we have seen so far.

The $\Delta$-activities can be defined in several manners. In this section, we are going to describe algorithms that compute such activities.

\subsection{Decision trees and decision functions}
\label{subsec:dectree}

All the families of activities we saw depended on a parameter: linear order for the ordering activities, embedding for the embedding activities... Similarly, the $\Delta$-activities will depend on a object, which is in a certain sense more general, named \emph{decision tree}.

\begin{figure}[h!]
\begin{center}
\begin{minipage}{3cm}\includegraphics[scale=1.15]{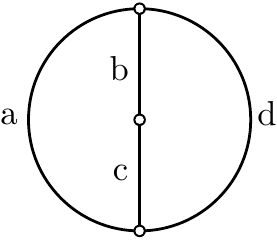}\end{minipage} \hspace{2cm} \begin{minipage}{7cm}\includegraphics[scale=1]{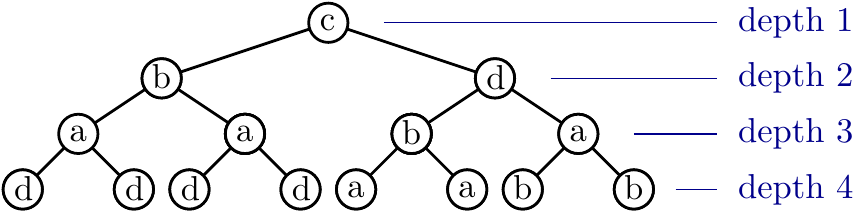}\end{minipage}
\end{center}
\caption{Example of a graph and an associated decision tree.}
\label{fig:ex}
\end{figure}

Consider a graph $G$. A \emph{decision tree} is a perfect binary tree\footnote{A \textit{perfect binary tree} is a binary tree in which every node has $0$ or $2$ children and all leaves are on the same level. Sometimes perfect trees are called full trees.} $\Delta$ with a labelling $V(\Delta) \rightarrow E(G)$
such that along every path starting from the root and ending on a leaf, the labels of the nodes form a permutation of $E(G)$. 
In particular, the depth of every leaf is $|E(G)|$. An example of decision tree is shown in Figure \ref{fig:ex}.

A \emph{direction} denotes either left or right, which we will write $\ell$ and $r$. Each node of the decision tree bijectively corresponds to a sequence $(d_1,\dots,d_k)$ of directions with $0 \leq k < |E(G)|$: the root node maps to the empty sequence, its children respectively map to the sequences $(\ell)$ and $(r)$, its grand-children map to $(\ell,\ell)$, $(\ell,r)$, $(r,\ell)$ and $(r,r)$, and so on. 
If a node $n$ maps to the sequence $(d_1,\dots,d_k)$, then the label of $n$ will be denoted by $\Delta(d_1,\dots,d_k)$. 
For instance, given the decision tree of Figure \ref{fig:ex}, we have $\Delta(r,\ell) = b$ and $\Delta(\ell,\ell,r) = d$.
By convention, the left or right child of a leaf in $\Delta$ is the empty node.

This  function $\Delta$ thus maps every sequence $(d_1,\dots,d_k)$ of directions to an edge of $G$ such that 
\begin{equation} \forall\, i < k \ \ \Delta(d_1,\dots,d_i) \neq \Delta(d_1,\dots,d_k),\end{equation}
for every sequence $(d_1,\dots,d_k)$ of directions. Such a function is called a
\emph{decision function}.

The decision function and the decision tree are the same object under different forms. Indeed, given a decision function, it is not difficult to label a decision tree accordingly. Therefore, when a formal definition of a decision tree is needed, we could give the decision function instead.

\subsection{Algorithm}

We fix $G$ a connected graph with $m$ edges and $\Delta$ a decision tree for $G$. Given a subgraph $S$ of $G$ (not necessarily a spanning tree), Algorithm \ref{type} outputs a partition ($TypeSe$, $TypeL$, $TypeSi$, $TypeI$) of $E(G)$. In other terms, the algorithm assigns to each edge a type, denoted by \bSe\, (for Standard External), \bL\, (for Loop), \bSi\, (for Standard Internal) or \bI\, (for Isthmus). Thus, if an edge belongs to $TypeSe$ (resp. $TypeL$, $TypeSi$, $TypeI$), we say that this edge has $\Delta$-\emph{type} \bSe\, (resp. \bL, \bSi, \bI) for $S$. If there is no ambiguity, we simply write \emph{type}. \\

\begin{algorithm}[h!]
\caption{How types are assigned to the edges.}
\label{type}
\begin{algorithmic}[5]
\Require $S$ subgraph of $G$.
\Ensure A partition ($TypeSe$, $TypeL$, $TypeSi$, $TypeI$) of $E(G)$.
\State $m \leftarrow $ number of edges in $G$;
\State $TypeSe \leftarrow \emptyset$; $TypeL \leftarrow \emptyset$; $TypeSi \leftarrow \emptyset$; $TypeI \leftarrow \emptyset$; 
\State $n \leftarrow$ root of $\Delta$;
\State $H \leftarrow G$;
\For {$k$ from $1$ to $m$}

	\State $e_{k} \leftarrow$ label of $n$;
	\If{$e_k$ is standard in $H$ \textbf{and} $e_k \notin S$ }

			\State $H \leftarrow \delete H e$;  
			\State \textbf{add} $e_k$ in $TypeSe$;
			\State $n \leftarrow$ left child of $n$;
	\EndIf
	\If{$e_k$ is a loop in $H$}
		\State {\color{darkgray} $H \leftarrow \delete H e$;} 			    
		\label{suppressionfacultative} 
		\Comment{optional (see Prop. \ref{variants})}
		\State \textbf{add }$e_k$ in $TypeL$;
		\State $n \leftarrow$ left child of $n$;
	\EndIf 
	\If{$e_k$ is standard in $H$ \textbf{and} $e_k \in S$}
		 	\State $H \leftarrow \contract H e$;
			\State \textbf{add} $e_k$ in $TypeSi$;
			\State $n \leftarrow$ right child of $n$;
	\EndIf 
	\If{$e_k$ is an isthmus in $H$}
		\State {\color{darkgray} $H \leftarrow \contract H e$ ;}
		  \label{contractionfacultative} 
		\Comment{optional (see Prop. \ref{variants})}
		\State \textbf{add }$e_k$ in $TypeI$;
		\State $n \leftarrow$ right child of $n$;
	\EndIf 
\EndFor
\State \Return ($TypeSe$, $TypeL$, $TypeSi$, $TypeI$)
\end{algorithmic}
\end{algorithm}

\noindent \textbf{Informal description.} We start from the edge that labels the root of $\Delta$. If this edge is standard external or a loop, we assign it the type \bSe\, or \bL, the edge is deleted and we go to the left subtree of $\Delta$. If this edge is standard internal or an isthmus, we assign the type \bSi\, or \bI, the edge is contracted and we go to the right subtree of $\Delta$.  We repeat the process until the graph has no more edge. Figure \ref{schema} illustrates this description.

\fig{[scale=1.4]{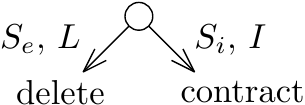}}{Diagram representing a step of Algorithm \ref{type}.}{schema}

We say that an edge is $\Delta$-\emph{active} (or \textit{active})  for $S$ if it has $\Delta$-type \bL\, or \bI. The $\Delta$-activity is the function denoted by $\Act$ that maps each spanning tree onto its set of $\Delta$-active edges.   Warning: an edge with type \bL\, or \bI\, is not necessarily a loop or an isthmus. It is an edge that is a loop or an isthmus at some point in Algorithm \ref{type}.  

\begin{figure}[h!]

\begin{center}
\hrule \  \\
$H$  at each iteration when   line \ref{suppressionfacultative} and line \ref{contractionfacultative} are present: \\
\vspace{0.1cm}

\includegraphics[scale=0.82]{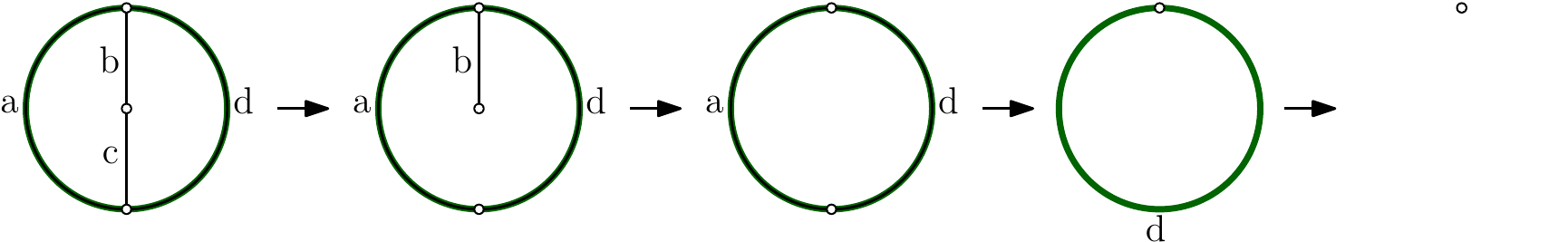} \\
\vspace{0.1cm}
\hrule \  \\
\vspace{0.1cm}

$H$ at each iteration when line \ref{suppressionfacultative} and line \ref{contractionfacultative} are missing: \\
\vspace{0.1cm}
\includegraphics[scale=0.82]{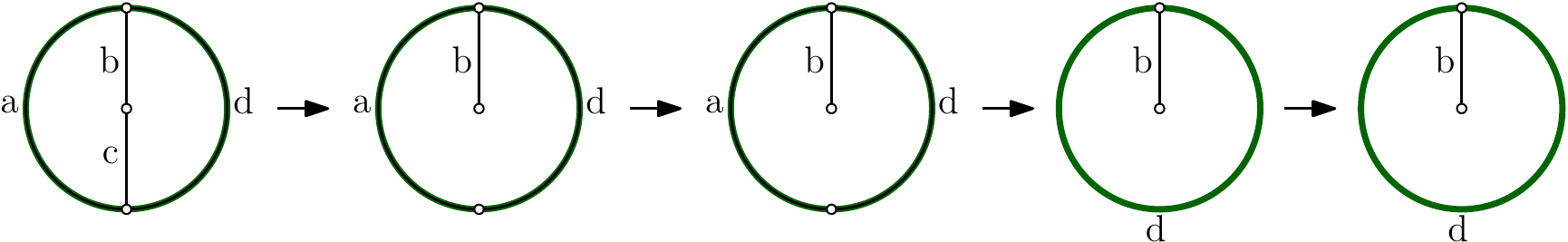}
\vspace{0.1cm}
\hrule
\end{center}

\caption{Run of two versions of Algorithm \ref{type} with $S = \{a,d\}$. One has assigned to $c$ the type \bSe, to $b$ the type \bI, to $a$ the type \bSi, to $d$ the type \bL, in this order. }
\label{ex:run}
\end{figure}

\noindent \textbf{Example.} Consider $G$ the graph and $\Delta$ the decision tree of Figure \ref{fig:ex}. The run of the algorithm for $S = \{a,d\}$ is illustrated by Figure \ref{ex:run} (top). We have $\Act(S) = \ens{b,d}$. 

Algorithm \ref{type} is a generalization of the matroid resolution algorithm from Gordon and Traldi  \cite{gordon-traldi} that computes the ordering activity. In Gordon and Traldi's version, the edges $(e_k)$ are considered in a fix order that does not depend on the subgraph. This constitutes a noticeable difference with Algorithm \ref{type} where the sequence $(e_k)$ is given by the decision tree.



A sequence of edges $(e_1,\dots,e_m)$ and a sequence of types $(t_1,\dots,t_m)$  can be naturally assigned to each subgraph $S$, where $m$ is the number of edges in $G$, where $e_k$ is the $k$-th edge visited by Algorithm \ref{type} and $t_k$ is the type of $e_k$. This pair of sequences is called the \emph{history} of $S$ and is denoted by $e_1 \ua{t_1} e_2 \ua{t_2} \dots  \ua{t_{m-1}} e_m \ua{t_{m}}$. For instance, given $G$ and $\Delta$ as in Figure \ref{fig:ex}, the history of $\{a,d\}$ is $c \ua{\mSe} b \ua{\bI} a \ua{\mSi} d \ua{\bL}$.

We can easily prove by induction that for every history $e_1 \ua{t_1} e_2 \ua{t_2} \dots  \ua{t_{m-1}} e_m \ua{t_{m}}$, the relation
\begin{equation}
\label{ekdelta}
e_{k+1} = \Delta(d(t_1),\dots,d(t_k))
\end{equation}
holds for every $k \in \ens{0,\dots,m-1}$, where $d$ denotes the map $\ens{\mSe,\bL,\mSi,\bI} \rightarrow \ens{\ell,r}$ defined as
\begin{equation} 
d(t) = \left\{\begin{array}{cl} \ell &   \textrm{if }t = \mSe\,\textrm{ or }t =\bL, \\ r &  \textrm{if }t = \mSi\,\textrm{ or }t =\bI. \end{array} \right.\end{equation}



\subsection{Description of the Tutte polynomial with $\Delta$-activity}


We are now ready to state the central theorem.

\begin{theo} \label{charact}
Let $G$ be any connected graph and $\Delta$ a decision tree for $G$. 
The $\Delta$-activity is Tutte-descriptive. In other terms, the Tutte polynomial of $G$ is equal to
\begin{equation}
\label{celuila}
T_G(x,y) = \sum_{T\textrm{ spanning tree of }G} x^{|\mathcal I(T)|}y^{|\mathcal E(T)|},
\end{equation}
where $\mathcal I(T)$ and $\mathcal E(T)$ are respectively the sets of internal $\Delta$-active and external $\Delta$-active edges of the spanning tree $T$.\end{theo}

Observe that this theorem holds for every decision tree $\Delta$, although the notion of $\Delta$-activity depends on the chosen decision tree.

Before giving a proof of this theorem, let us remark that when the subgraph is a spanning tree, the edges of type \bI/\bL\, coincide with the internal/external active edges.
\begin{prop}
\label{intextact}
Let $G$ be a connected graph and $\Delta$ a decision tree for $G$. 
If the considered subgraph $S=T$ is a spanning tree of $G$, then each edge of type \bI\, is internal and each edge of type \bL\, is external.
\end{prop}
\begin{proof}
It is not difficult to see that at each step of the algorithm, $T \cap H$ is a spanning tree of $H$. Thus, when the algorithm assigns to an edge $e$ the type \bI\, (resp. \bL), $e$ is an isthmus (resp. a loop) in $H$, so it must be inside $T$ (resp. outside $T$).
\end{proof}

 In other terms, for every spanning tree $T$, the edges of type \bI\, are precisely the internal $\Delta$-active edges  and the edges of type \bL\, the external $\Delta$-active edges. Also note that each edge of type \bSe\, is external and each edge of type \bSi\, is internal, even if the subgraph is not a spanning tree.
 

\begin{proof}[Proof of Theorem \ref{charact}] Let $\tut$ be the polynomial 
$$\tut(x,y)=\sum_ {T\textrm{ spanning tree of }G} x^{|\mathcal I(T)|}y^{|\mathcal E(T)|}.$$
Let us prove by induction on the number of edges that the equality $T_G = \tut$ is true for any graph $G$ and decision tree $\Delta$.
More precisely, we show that $\tut$ satisfies the same relation of induction as the Tutte polynomial (see Proposition~\ref{eq:ind}).

 Assume that $G$ is reduced to a graph with one edge $e$. Then $e$ is a loop or an isthmus, and $\Delta$ is reduced to a leaf labelled by $e$. We easily check that    $T_G(x,y)=y=\tut(x,y)$ when $e$ is a loop and $T_G(x,y)=x=\tut(x,y)$ when $e$ is an isthmus.

 Now assume that $G$ has at least two edges. We denote by $n_1$ the root node of $\Delta$ and $e_1$ its label.

\begin{lem} \label{respect}
Let $T$ be a spanning tree of $G$. If $e_1$ is external for $T$, define the graph $G'$ as $\delete G {e_1}$
and $\Delta'$ as the subtree rooted on the left child of $n_1$. If $e_1$ is internal, define the graph $G'$ as $\contract G {e_1}$ and $\Delta'$ as the subtree rooted on the right child of $n_1$.
Then the $\Delta$-type in $G$ for $T$ and the $\Delta'$-type in $G'$ for $T \backslash e_1$ of every edge in $G'$ are the same. 
\end{lem}

\begin{proof} 
After the first iteration of Algorithm \ref{type} with graph $G$, decision tree $\Delta$ and input $T$, one can check that the graph $H$ is equal to $G'$, and the node $n$ is the root node of $\Delta'$. (Indeed, by Proposition \ref{intextact}, an edge in a spanning tree is external if and only if it has type \bSe\, or \bL. So we go to the left subtree of $\Delta$ at the first iteration if and only if $e_1$ is external.) These are exactly the values of $H$ and $n$ at the beginning of Algorithm \ref{type} with graph $G'$, decision tree $\Delta'$ and input $T \backslash e_1$. Therefore, in both cases, the algorithm will evolve in the same way from this point on: the edges will be visited in the same order and the types will be identically assigned.
%
\end{proof}

We study now three  cases:

\noindent \textbf{(1) $\boldsymbol{e_1}$ is a standard edge.} Then we can partition the set of spanning trees of $G$ into two disjoint sets: the set of spanning trees \textit{not} containing the edge $e_1$, denoted by $\mathbb T_1$, and the set of spanning trees  containing $e_1$, denoted by $\mathbb T_2$. Let us consider $\Delta_1$ (resp. $\Delta_2$) the subtree rooted on the left child (resp. right child) of $n_1$.

The map $\phi : T \mapsto \delete T {e_1}$ is a bijection from $\mathbb T_1$ to the spanning trees of $\delete G  {e_1}$. (Actually $\phi$ is nothing else than the identity map.) Moreover, for every $T \in \mathbb T_1$, the edge $e_1$ has type \bSe. So according to Lemma \ref{respect}, the bijection $\phi$ preserves the internal and external active edges (for $G$ and $\Delta$ on the one hand, for $G'$ and $\Delta_1$ on the other hand), hence
$$\mathcal T_{\delete G {e_1},\Delta_1} = \sum_{T \in \mathbb T_1} x^{|\mathcal I(T)|}y^{|\mathcal E(T)|}.$$
Similarly, one can prove that 
$$\mathcal T_{\contract G {e_1},\Delta_2} = \sum_{T \in \mathbb T_2} x^{|\mathcal I(T)|}y^{|\mathcal E(T)|}.$$
But we have $\mathcal T_{\delete G {e_1},\Delta_1} = T_{\delete G {e_1}}$ and $\mathcal T_{\contract G {e_1},\Delta_2} = T_{\contract G {e_1}}$ by the induction hypothesis.
So with Proposition \ref{eq:ind}, we get
\begin{eqnarray*} T_G(x,y) & = & T_{\delete G {e_1}}(x,y) + T_{\contract G {e_1}}(x,y) \\ \ & = & \sum_{T \in \mathbb T_1} x^{|\mathcal I(T)|}y^{|\mathcal E(T)|} + \sum_{T \in \mathbb T_2} x^{|\mathcal I(T)|}y^{|\mathcal E(T)|} = \tut(x,y).
\end{eqnarray*}

\noindent \textbf{(2) $\boldsymbol{e_1}$ is an isthmus.} Then $e_1$ must be internal for every spanning tree of $G$ and so the map $T \mapsto \contract T {e_1}$ is a bijection from the spanning trees of $G$ to the spanning trees of $\contract G {e_1}$. 
Moreover, for each spanning tree $T$ of $G$, the edge $e_1$ has $\Delta$-type \bI\, because it is an isthmus. Let $\Delta'$ denote the subtree of $\Delta$ rooted on the right child of $n_1$. By Lemma \ref{respect}, the $\Delta$-type for $T$ and the $\Delta'$-type for $\contract T {e_1}$ are identical  for every edge different from $e_1$.
So from the properties above, we have $$\tut(x,y) = x\, \mathcal T_{\contract G {e_1},\Delta'}(x,y).$$ We conclude by induction and Proposition \ref{eq:ind}.

\noindent \textbf{(3) $\boldsymbol{e_1}$ is a loop.} Exactly the dual demonstration of the case (2).
\end{proof}

\subsection{Variants of the algorithm}

Algorithm \ref{type} is rather flexible. As stated inside the pseudo-code, it admits variants that lead to the same partition of edges.

\begin{prop} \label{variants}
Removing line \ref{suppressionfacultative} or/and line  \ref{contractionfacultative} from Algorithm \ref{type} does not change its output. In other terms, given a subgraph $S$ of $G$, the types of the edges remain identical if we choose not to delete the edges of type \bL\, or/and not to contract the edges of type \bI.
\end{prop}

\begin{proof} When a graph $G'$ is obtained by deletions and contractions of some edges from a graph $G$, we write $G \rightsquigarrow G'$. In this case, note that if $e$ is a loop (resp. an isthmus) in $G$ and an edge in $G'$, then $e$ is still a loop (resp. an isthmus) for $G'$.
Let $A$ be the variant of Algorithm \ref{type} where the lines   \ref{suppressionfacultative} and \ref{contractionfacultative} are missing and  $A'$ a variant where these lines are potentially present.
The graphs $H$ in the algorithms $A$ and  $A'$ just before the $j$-th iteration will be respectively  denoted by  $H_j$ and  $H'_j$.
Given a subgraph $S$ of $G$, let $e_1 \ua{t_1}  \dots  \ua{t_{m-1}} e_m \ua{t_{m}}$ (resp. $e'_1 \ua{t'_1}  \dots  \ua{t'_{m-1}} e'_m \ua{t'_{m}}$) be the history of $S$ for the algorithm $A$ (resp. $A'$). Let us prove by induction on $k$ that for each $1 \leq j \leq k$, we have $e_j = e'_j$ and $t_j = t'_j$. 

Assume that the induction hypothesis is true for $k-1$. The edges $e_k$ and $e'_k$ are identical since they are equal to $\Delta(d(t_1),\dots,d(t_{k-1}))$ (see Equation \eqref{ekdelta}). 
Let us prove the equivalence 
$$e_k\textrm{ loop in } H_k \Leftrightarrow e_k\textrm{ loop in } H'_k.$$
 The left-to-right implication is obvious since $H_k \rightsquigarrow H'_k$.
Conversely, if $e_k$ is not a loop in $H_k$, then there exists a cocycle $C$ in $H_k$ including $e_k$. Let us show that for every edge $x$ of $H_k$, we have $x \notin H'_k \Rightarrow x \notin C$, which implies that $C$ is included in $H'_k$ and so $e_k$ is not a loop in $H'_k$. Let $e$ be an edge of $H_k \backslash H'_k$. It must be an edge of the form $e_j$ with $j < k$ which has type \bI\, or \bL. So this edge was an isthmus or a loop in $H_j$, hence also in $H_k$ since $H_j \rightsquigarrow H_k$. Therefore $e$ cannot appear in a non-singleton cocycle of $H_k$, and in particular in $C$.


Similarly we prove 
$$e_k\textrm{ isthmus in }H_k \Leftrightarrow e_k\textrm{ isthmus in }H'_k.$$

So if $e_k$ is respectively a loop, an isthmus, a standard external edge, a standard internal edge in  $H_k$, then it will be a loop, an isthmus, a standard external edge, a standard internal edge in $H'_k$; hence $t_k = t'_k$.\end{proof}

Each version of Algorithm \ref{type} can be of interest: For implementation, it would be better to perform a minimum number of operations and consequently choose the variant where the edges of type \bL\, or \bI\, remain untouched. From a theoretical point of view, deleting each edge of type \bL\, and contracting each edge of type \bI\, can facilitate the proofs. For example, with this version, it is easy to see that $e_m$ has necessarily type \bL\, or \bI\, since at the last iteration $e_m$ is the only edge of $H$, so must be an isthmus or a loop.

Moreover, Algorithm \ref{type} can be differently  adapted depending on the context. For instance, Algorithm \ref{algint} allows us to compute the internal $\Delta$-active edges of a spanning tree, the external active edges being omitted. This algorithm will be useful to connect blossoming activities with $\Delta$-activities.

\begin{algorithm}[h!]
\caption{Another way to compute the set of internal active edges.}
\label{algint}
\begin{algorithmic}[5]
\Require $T$ spanning tree of $G$.
\Ensure A subset $IntAct$ of $E(G)$.
\State $m \leftarrow $ number of edges in $G$; $IntAct \leftarrow \emptyset$; 
\State $n \leftarrow$ root of $\Delta$; $H' \leftarrow G$
\For {$k$ from $1$ to $m$}

	\State $e_{k} \leftarrow$ label of $n$;
	\If{$e_k \notin T$ and $e_k$ is not an isthmus in $H'$}
			\State $H' \leftarrow \delete {H'} e$;  
			\State $n \leftarrow$ left child of $n$;
	\Else
   			\State $n \leftarrow$ right child of $n$;
	\EndIf

	\If{$e_k$ is an isthmus in $H'$}
		\State \textbf{add }$e_k$ in $IntAct$;
	\EndIf 
\EndFor
\State \Return $IntAct$
\end{algorithmic}
\end{algorithm}

\noindent \textbf{Informal description:} The principle is the same as in Algorithm \ref{type}, except that the internal edges are never contracted.

\begin{prop} \label{prop:algint}
For any spanning tree $T$ of $G$, Algorithm \ref{algint} outputs the set of edges of type \bI.
\end{prop}
\begin{proof}
On the one hand, we consider the version of Algorithm \ref{type} where the edges of type \bI\, are not contracted and where the edges of type $\bL$ are deleted. Thus, the graph $H$ at the $k$-th iteration is obtained by contracting all the edges $e_i$ of type \bSi\, and by deleting all the external edges $e_i$, with $i \in \ens{1,\dots,k-1}$\footnote{By Proposition \ref{intextact}, an edge is external  in a spanning tree if and only if it has type \bSe\, or \bL.}. On the other hand, in Algorithm \ref{algint}, the graph $H'$ at the $k$-th iteration is obtained by deleting all the external edges of the form $e_i$, with $i \in \ens{1,\dots,k-1}$. So $H'$ differs from $H$ at this moment only by some edge contractions (which do not involve $e_k$). Therefore, $e_k$ is an isthmus in $H'$ if and only if $e_k$ is an isthmus in $H$. This means that $e_k$ belongs to the output of Algorithm \ref{algint} if and only if $e_k$ has type \bI.
\end{proof}

\section{Edge ordering and $\Delta$-activity}
\label{s:ord}

It turns out that we also can define the notion of $\Delta$-activity by using fundamental cycles and cocycles, as Tutte and Bernardi did. This is the purpose of this section.

\subsection{$(\Delta,S)$-ordering}

Consider a graph $G$ with decision tree $\Delta$. Let $e_1 \ua{t_1} e_2 \ua{t_2} \dots  \ua{t_{m-1}} e_m \ua{t_{m}}$ denote the history of a subgraph $S$ of $G$. We define \textit{the $(\Delta,S)$-ordering} on the edge set of $G$ by setting:
$$e_1 < e_2 < \dots < e_m.$$
It corresponds to the visit order of the edges in Algorithm \ref{type}. This notion will be especially used for the subgraphs $S=T$ that are spanning trees.

Given a decision tree $\Delta$ and a spanning tree $T$ of a graph $G$, the $(\Delta,T)$-ordering can be easily constructed from $\Delta$ and $T$: Start with the root node of $\Delta$. If the edge label is external, go  down to the left child. If the edge is internal, go down to the right child. Repeat this operation until joining a leaf node. The sequence of the node labels gives the $(\Delta,T)$-ordering.

The explanation is simple. It relies on 
$\eqref{ekdelta}$ and the fact that for all spanning trees, an edge is external if and only it has type \bSe\, or \bL (see Proposition \ref{intextact}).

\noindent \textbf{Example.} Consider the graph and the decision tree $\Delta$ from Figure \ref{fig:ex}, with the spanning tree $T=\ens{a,c}$. By following the path corresponding to the sequence of directions $(r,\ell,\ell)$, we can see that the $(\Delta,T)$-ordering is equal to $c < d < b < a$.

\subsection{Fundamental cycles and cocycles}

The following proposition is the first step to prove the link between the $\Delta$-activity and some activities from Section \ref{sec:activity}.
\begin{prop} 
\label{maximal}
Consider a graph $G$ with decision tree $\Delta$ and a spanning tree $T$. An external (resp. internal) edge $e$ is $\Delta$-active if and only if it is maximal for the $(\Delta,T)$-ordering in its fundamental cycle (resp. cocycle).
\end{prop}

\noindent \textbf{Remark.} The active edges are here characterized by maximality in their corresponding fundamental cycles/cocycles, and not by minimality as in Tutte's or Bernardi's works. We cannot simply reverse the order so that  "maximal" becomes "minimal". Indeed, we would need the reverse order to correspond to a $(\Delta,T)$-ordering, which seldom happens (except for the ordering activity). This is related to the notion of tree-compatibility, which is treated in the next subsection.

Before proving this proposition, let us state a lemma that will be used on several times in this paper. 
\begin{lem}
\label{cyclecocycle}
For each subgraph $S$ of $G$, we have the following properties:
\begin{enumerate}
	\item For each edge $e$ of type \bL, there exists a cycle in $G$ which consists of $e$ and edges of type \bSi. Moreover, $e$ is maximal in this cycle for the $(\Delta,S)$-ordering.
	\item For each edge $e$ of type \bI, there exists a cocycle of $G$ which consists of $e$ and edges of type \bSe. Moreover, $e$ is maximal in this cocycle for the $(\Delta,S)$-ordering.
\end{enumerate}
\end{lem}
\begin{proof} Here we consider the variant of Algorithm \ref{type} where the edges of type \bI\, are \textit{not} contracted and the edges of type \bL\, are \textit{not} deleted.

(1) Suppose that $e=e_k$, i.e. $e$ was visited at the $k$-th iteration of the algorithm. Then let us prove by a decreasing induction on $j \in \{k,\dots,1\}$ that the graph $H$ at the $j$-th iteration contains a path linking the endpoints of $e$ only made of edges of type \bSi\, and visited before $e$. 
For the base case $j=k$, $e_k$ is a loop, so the empty path matches. Then assume that the induction hypothesis is true for some $j \in \{k,\dots,2\}$. The only operation which could ensure that the endpoints of $e$ can be linked with edges of type \bSi\, after the $(j-1)$-th iteration, but not before, is edge contraction. But a contracted edge must have type \bSi. So, whatever the type of $e_{j-1}$ is, there will be always, at the $(j-1)$-th iteration, a path made of edges of type \bSi\, and visited before $e$ linking the endpoints of $e$. So the induction step holds: the resulting path with the edge $e$ forms the expected cycle.

(2) Similar to point (1), in a dual way.
\end{proof}

We can now tackle the proof of Proposition \ref{maximal}.

\begin{proof}[Proof of Proposition \ref{maximal}] 1. Assume that an external edge $e_j$ is maximal in its fundamental cycle $C$. We  use a version of Algorithm \ref{type} where edges of type \bI\, are contracted. Except $e_j$, the fundamental cycle $C$ of $e_j$ is made of internal edges, so edges with type \bSi\, or \bI. Thus, at the $j$-th iteration, the edges other than $e$ have been contracted, which implies that $e_j$ is a loop. Therefore $e_j$ has type $\bL$. In other terms, it is an external active edge.

2. Assume that an external edge $e_j$ is active, that is to say it has type \bL. By Lemma \ref{cyclecocycle}, there exists a cycle $C$ made of $e$ and edges of type \bSi. Since the edges of type \bSi\, are internal, $C$ is the fundamental cycle of $e$. Lemma \ref{cyclecocycle} also claims that $e$ has been visited last in $C$: it means that $e$ is maximal in $C$ for the $(\Delta,T)$-ordering.

The case where $e$ is internal can be processed in the same way.
\end{proof}

\subsection{Order map and tree-compatibility}

If we want to prove that a specific activity (for example an embedding activity) is a $\Delta$-activity, we need to define the adequate decision tree. However, building this decision tree can be rather tricky. In this subsection, we are going to give a combinatorial condition that guarantees the existence of such a decision tree. 

Fix a graph $G$. An \emph{order map} is a map from the set of spanning trees of $G$ onto the set of total orders of $E(G)$ (or equivalently the permutations of $E(G)$). Given an order map $\phi$ and a spanning tree $T$, we denote by $\phi_k(T)$ the $k$-th smallest edge in this order.

An order map $\phi$ is said to be \emph{tree-compatible} if we can construct a decision tree $\Delta_\phi$ such that for every spanning tree $T$, the $(\Delta_\phi,T)$-ordering corresponds to $\phi(T)$. In other, a tree-compatible order map is such that for each $k \in \ens{1,\dots,|E(G)|}$, the edge $e_k$, defined by Algorithm \ref{type} with input $T$ for some decision tree $\Delta$, is equal to $\phi_k(T)$. 

\noindent \textbf{Example.} Consider the graph from Figure \ref{fig:ex} with this order map :
$$ \phi : \begin{array}{ccc}
\ens{a,b} & \mapsto & c < b < a < d \\
\ens{a,c} & \mapsto & c < d < b < a \\
\ens{b,c} & \mapsto & c < d < b < a \\
\ens{b,d} & \mapsto & c < b < a < d \\
\ens{c,d} & \mapsto & c < d < a < b
\end{array} .
$$
The order map $\phi$ is tree-compatible: if we denote by $\Delta_\phi$  the decision tree of Figure \ref{fig:ex}, then we can easily check that $\phi(T)$ and the
$(\Delta_\phi,T)$-ordering are the same for all spanning trees $T$.
We can notice that the decision tree is not unique. For example, the same decision tree fits if we replace the leftmost subtree with three nodes, namely \begin{tabular}{l}\includegraphics[height=20px]{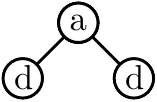}\end{tabular}, by \begin{tabular}{l}\includegraphics[height=20px]{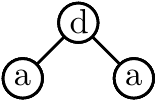}\end{tabular}.


The following theorem gives a characterization of tree-compatible order maps.

\begin{theo} \label{treecompatible}
An order map $\phi$ is tree-compatible if and only if 
for all spanning trees $T$ and $T'$ and for every $k \in \ens{0,\dots,|E(G)|-1}$ the implication 
\begin{multline} \label{tc}
T \cap \ens{\phi_1(T),\dots,\phi_k(T)} = T' \cap \ens{\phi_1(T),\dots,\phi_k(T)} \\ \Longrightarrow \ \forall j \in \ens{1,\dots,k+1} \ \, \phi_j(T) = \phi_j(T')
\end{multline} holds.

\end{theo}

Note the case $k = 0$ is significant: it implies that $\phi_1(T)$ is the same edge for all spanning trees $T$.

\noindent \textbf{Back to the previous example.} Consider the order map $\phi$ defined by \eqref{exom} and set $T = \ens{a,c} $ and $T' = \ens{c,d}$. We have $\phi_1(T) = c$ and $c$ belongs to $T$ and $T'$. So for $k = 1$, the antecedent of \eqref{tc} holds. Thus the order map coincides one step further: \mbox{$\phi_2(T) = \phi_2(T')$} ($=d$). But for $k=2$, the antecedent of \eqref{tc} does not hold since $\phi_2(T) = d \notin T$ and $d \in T$: no inferences can be made from this observation.

\begin{proof} \textbf{Left-to-right implication.} Let $\Delta$ be a decision tree and $\phi$ be the order map that sends each spanning tree $T$ onto the $(\Delta,T)$-ordering. Let us prove by induction on $k$ that \eqref{tc} is true for all spanning trees $T$ and $T'$. It holds for $k = 0$ since $\phi_1(T)$ is for all spanning trees $T$ the label of the root node of $\Delta$ (which is the first edge we visit in Algorithm \ref{type}). Assume now that 
$$T \cap \ens{\phi_1(T),\dots,\phi_{k+1}(T)} = T' \cap \ens{\phi_1(T),\dots,\phi_{k+1}(T)}.$$ By the induction hypothesis,  \eqref{tc} holds, hence for all $j \in \ens{1,\dots,k+1}$ we have $\phi_j(T) = \phi_j(T')$. Moreover, by  \eqref{ekdelta}, we have $\phi_{k+2}(T) = \Delta(d^T_1,\dots,d^T_{k+1})$, where $d^T_i = \ell$ when $\phi_i(T)$ is external in $T$ and $d^T_i = r$ when $\phi_i(T)$ is internal in $T$. (Recall that the type an edge $e_i$ in a spanning tree is equal to \bSe\, or \bL\, when $e$ is external, and is equal to \bSi\, or \bI\, when $e$ is internal.) Of course, it also holds if we replace $T$ by $T'$. But by assumption, for every $i \in \, \ens{1,\dots,k+1}$, the edge $\phi_i(T)$ is internal in $T$ if and only if $\phi(T')$ is internal in $T'$. Consequently, we have $\phi_{k+2}(T) = \phi_{k+2}(T')$.

\textbf{Right-to-left implication.} Let $\phi$  be an order map that satisfies \eqref{tc} for  all spanning trees $T$ and $T'$ and for every $k \in \ens{0,\dots,|E(G)|-1}$.
As indicated in Subsection \ref{subsec:dectree}, we can define a decision function instead of a decision tree in order to prove that $\phi$ is tree-compatible. 

\textbf{(1) Inductive definition of the decision function.} Consider a sequence of directions $(d_1,\dots,d_{k-1})$ with $k \in \ens{1,\dots,|E(G)|}$. We suppose by induction that for $i \in \ens{1,\dots,k-1}$ the edge $\Delta_\phi(d_1,\dots,d_{i-1})$, denoted by $\eta_{i}$, is well defined. Given a spanning tree $T$ and  $k \in \ens{1,\dots,|E(G)|}$, we denote by $P(T,k)$ the property
\begin{equation*} 
\ens{ j \in \ens{1,\dots,k-1} \  |  \  \eta_j \in T } = \ens{ j \in \ens{1,\dots,k-1} \  |  \ d_j = r }.
\end{equation*}
Always by induction, we assume that for every $i \in \ens{1,\dots,k-1}$
  if there exists a spanning tree $T$ that satisfies $P(T,i)$, then $\eta_i$ equals $\phi_i(T)$.
  (The base case of the induction is embodied by the case $k=1$. Indeed, when $k=1$, the set $\ens{1,\dots,k-1}$ is empty : no induction hypothesis is assumed.) 

If there exists a spanning tree $T$ that satisfies $P(T,k)$, we define $\Delta_\phi(d_1,\dots,d_{k-1})$ as $\phi_k(T)$. Two points have to be checked :
\begin{enumerate}
\item[(a)] $\Delta_\phi(d_1,\dots,d_{k-1})$ is  different from $\eta_1$, $\dots$, $\eta_{k-1}$. For every $i \in \ens{1,\dots,k-1}$, the property $P(T,i)$ holds because $P(T,k)$ holds. So by the induction hypothesis, we have $\eta_i = \phi_i(T)$, which is well different from $\Delta_\phi(d_1,\dots,d_{k-1}) = \phi_k(T)$.

\item[(b)] The definition of $\Delta_\phi(d_1,\dots,d_k)$ does not depend on the chosen spanning tree. Suppose that there exist two spanning trees $T$ and $T'$ such that $P(T,k)$ and $P(T',k)$ hold. Then $\ens{ i \in \ens{1,\dots,k-1} \  |  \  \eta_i \in T } = \ens{ i \in \ens{1,\dots,k-1} \  |  \  \eta_i \in T'}$. By the induction hypothesis,e have $\eta_i = \phi_i(T)$ for all $i \in \ens{1,\dots,k-1}$, hence $T \cap \ens{\phi_1(T),\dots,\phi_{k-1}(T)} = T' \cap \ens{\phi_1(T),\dots,\phi_{k-1}(T)}$. By \eqref{tc}, we have $\phi_k(T) = \phi_k(T')$. 
\end{enumerate}

If there exists no spanning tree $T$ such that $P(T,k)$ is true, we arbitrarily define $\Delta_\phi(d_1,\dots,d_{k-1})$ as an edge different from $\eta_1$, $\dots$, $\eta_{k-1}$. Here there is no importance for the choice of the edge, since the corresponding branch of the decision tree will never  be  visited.



\textbf{(2) The $\boldsymbol{(\Delta_\phi,T)}$-ordering coincides with $\boldsymbol \phi$.} Let $k \in \ens{1,\dots,|E(G)|}$ and $T$ be a spanning tree. By  \eqref{ekdelta}, the edge $e_k$ (set by Algorithm \ref{type} with decision tree $\Delta_\phi$ and input $T$) equals $\Delta_\phi(d_1,\dots,d_{k-1})$, the direction $d_i$ being defined for $i \in \ens{1,\dots,k-1}$ as:
\begin{equation} \label{defdi}
d_i = \left\{\begin{array}{cl} \ell &   \textrm{if }e_i\textrm{ is external for }T, \\ r & \textrm{if }e_i\textrm{ is internal for }T. \end{array} \right.\end{equation}
Hence,
\begin{align*}
\ens{ i \,  \in \{1,\dots,k-1\} \  |  \  \Delta_\phi(d_1,\dots,d_{i-1}) \in T } & = \ens{i  \,  \in \{1,\dots,k-1\} \  |  \ e_i \in T } \\
\ & = \ens{i \,  \in \{1,\dots,k-1\} \  |  \ d_i = r }.
\end{align*}
So the spanning tree $T$ satisfies $P(T,k)$. By definition of $\Delta_\phi$, we have $$e_k = \Delta_\phi(d_1,\dots,d_{k-1}) = \phi_k(T),$$ which proves the tree-compatibility of $\phi$. \end{proof}


The following corollary will be helpful to prove that the activities from Section~\ref{sec:activity}  are Tutte-descriptive. Observe that its statement does not involve Algorithm \ref{type}.

\begin{cor} \label{cor:act} Let $\phi$ be an order map. Assume that $\phi$ is tree-compatible, i.e. for all spanning trees $T$ and $T'$ and for $k \in \ens{0,\dots,|E(G)|-1}$ such that $T \cap \ens{\phi_1(T),\dots,\phi_k(T)} = T'~\cap~\ens{\phi_1(T),\dots,\phi_k(T)}$, we have $\phi_i(T)=\phi_i(T')$ for each $i \in \ens{1,\dots,k+1}$. 

The activity that maps any spanning tree $T$ onto the set of edges that are maximal for $\phi(T)$ in their fundamental cocycles/cycles is a $\Delta$-activity. Therefore it is Tutte-descriptive.
\end{cor}
\begin{proof} By the definition of tree compatibility, there exists a decision tree $\Delta_\phi$ such that  the $(\Delta_\phi,T)$-ordering is  $\phi(T)$ for every spanning tree $T$. Furthermore, Proposition \ref{maximal} states that an edge is $\Delta_\phi$-active if and only if it is maximal in its fundamental cycle/cocyle for the $(\Delta_\phi,T)$-ordering, that is, $\phi(T)$. Then we conclude thanks to Theorem~\ref{charact}.
\end{proof}

\section{Partition of the subgraph poset into intervals }
\label{s:partition}

Crapo discovered in \cite{crapo} that the ordering activity induces a natural partition of the poset of the subgraphs with some nice properties. Bernardi, Gessel and Sagan also defined a similar partition based on their respective notions of activity. In this section, we prove the universality of this partition by extending it to the $\Delta$-activities.

\subsection{Introduction on an example}

We are going to introduce the different properties of this section throughout an example. Consider $G$ the graph and $\Delta$ the decision tree of Figure \ref{fig:ex}. The following table lists the types of edges for all subgraphs of $G$ (described by their sets of edges).

\begin{center}
\noindent \begin{tabular}{|c|c|c|c|c|}
\hline Subgraphs & Type of $a$ & Type of $b$ & Type of $c$ & Type of $d$ 
\\
\hline

$\emptyset$, $\{b\}$, $\{d\}$, $\mathbf{\{b,d\}}$  & \bSi & \bI & \bSe & \bI \\
\hline
$\{a\}$, $\mathbf{\{a,b\}}$, $\{a,d\}$, $\{a,b,d\}$  & \bSi & \bI & \bSe & \bL \\
\hline
$\{c\}$, $\mathbf{\{a,c\}}$  & \bI & \bSe & \bSi & \bSe \\
\hline
$\mathbf{\{b,c\}}$, $\{a,b,c\}$ & \bL & \bSi & \bSi & \bSe \\
\hline
$\mathbf{\{c,d\}}$, $\{a,c,d\}$, $\{b,c,d\}$, $\{a,b,c,d\}$ & \bL & \bL & \bSi & \bSi \\
\hline
\end{tabular}
\end{center}
\vspace{0.1cm}

\begin{figure}[h!]
\begin{center}
\includegraphics[scale=1]{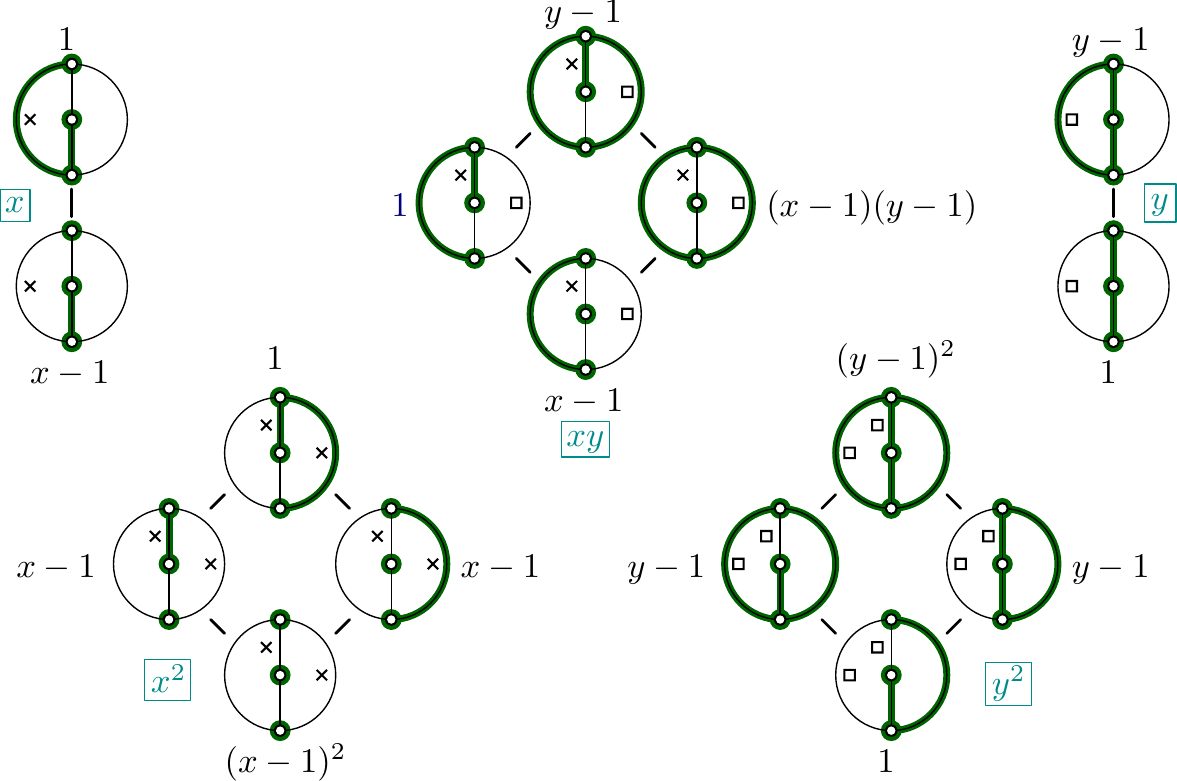}
\end{center}
\caption{Partition of the subgraphs with respect to the relation $\sim$, given the graph and decision tree from Figure \ref{fig:ex}. The crosses indicate the edges with type \bI\, and the squares the edges with type \bL. The contribution of each subgraph to the Tutte polynomial is also mentioned, as well as the contribution of each part of the partition.}
\label{fig:part}
\end{figure}

We have gathered in each line the subgraphs that share the same partition of edges. We can observe that the set of subgraphs inside any line is a subgraph interval. More particularly, given $S$ and $S'$ two subgraphs taken from the same line, $\diffs S S'$ is only made of active edges (for $S$ \emph{and} for $S'$). This will be proved in Proposition \ref{samehistory}. Furthermore, in each line we can find exactly one spanning tree, indicated in bold. This means that a spanning tree is included in each part of the partition of subgraphs, as stated in Theorem \ref{partition}. Finally, observe that if we put a weight $x$ per edge of type \bI\, and a weight $y$ per edge of type \bL, then by summing over all the lines of the table we get
$$x^2 + x\,y + x + y + y^2,$$
which exactly corresponds to the Tutte polynomial of $G$.  This is a consequence of Proposition \ref{rouge}. All these properties can be observed on Figure  \ref{fig:part}.

\subsection{Equivalence relation}


For the rest of the section, we fix $G$ a graph and $\Delta$ a decision tree. Given two subgraphs $S$ and $S'$ of $G$, we say that $S$ and $S'$ are \emph{equivalent} if they share the same history. In this case, we write $S \sim S'$. 
This is obviously an equivalence relation. Here are some other characterizations of this relation.

\begin{prop}
\label{samehistory}
Let $S$ and $S'$ be two subgraphs of $G$. Denote by $\Act(S)$ the set of active edges of $S$. The following properties are equivalent:
\begin{enumerate}
\item [(i)] $S \sim S'$,
\item [(ii)] the subgraphs $S$ and $S'$ induce the same partition of $E(G)$,
\item [(iii)] the subgraphs $S$ and $S'$ yield the same set of edges with type \bSe\, and the same set of edges with type \bSi,
\item [(iv)] $\diffs S {S'} \subseteq \Act(S)$,
\item [(v)] there exists $R \subseteq \Act(S)$ such that $S' = \diffs S R$.
\end{enumerate}


\end{prop}
\begin{proof} \textbf{(i) $\mathbf{ \Rightarrow }$ (ii) $\mathbf{ \Rightarrow }$ (iii)}. Trivial. \\
\noindent \textbf{(iii) $\mathbf{ \Rightarrow }$ (iv)}. Let the set of edges with type \bSi\, be denoted by $J$. So we have $S = J \cup (S \cap \Act(S))$ and $S' = J \cup (S' \cap \Act(S'))$ since an edge of type \bSi\, is always internal and an edge of type \bSe\, always external. Thus $$\diffs S {S'} = \diffs  {(S \cap \Act(S))} {(S' \cap \Act(S'))} \subseteq \Act(S) \cup  \Act(S').$$
But an edge that has neither type \bSe\, nor \bSi\, is active, so $\Act(S) = \Act(S')$ and then $\diffs S {S'} \subseteq \Act(S).$ \\
\textbf{(iv) $\mathbf{ \Rightarrow }$ (v)}. Set $R = \diffs S {S'}$.
Since the symmetric difference is associative, we have
$$\diffs S R = \diffs S {(\diffs S {S'})} = \diffs {(\diffs S S)} {S'} = S'.$$
\textbf{(v) $\boldsymbol{ \Rightarrow }$ (i)}: Let $e_1 \ua{t_1}  \cdots  \ua{t_{m-1}} e_m \ua{t_{m}}$ (resp. $e'_1 \ua{t'_1}  \cdots  \ua{t'_{m-1}} e'_m \ua{t'_{m}}$) be the history of $S$ (resp. $S'=\diffs S R$); we are going to prove by induction on $k$ that for each $1 \leq j \leq k$, we have $e_j = e'_j$ and $t_j = t'_j$. 
We assume that the inductive hypothesis is true at the step $k-1$. Then, by \eqref{ekdelta}, we have $$e_k = \Delta(d(t_1),\dots,d(t_{k-1})) = e'_k.$$
Let us consider the beginning of the $k$-th iteration.  The graph $H$ has undergone the same transformations whether the input was $S$ or $\diffs S R$. Indeed, the deletions or contractions we have performed are governed by the types of the visited edges, which are identical by the inductive hypothesis. 
If the edge $e_k$ is a loop (resp. an isthmus) in $H$, then by following the algorithm we get $t_k = t'_k = \bL\, $ (resp. $t_k = t'_k =\bI$). If $e_k$ is standard, then $e_k$ is not active for $S$ and so $e_k \notin R$. In this case, either  $e_k \notin S$ which implies $e_k \notin \diffs S R$ and $t_k=t'_k=\mSe$, or $e_k \in S$ which implies $e_k \in \diffs S R$ and $t_k=t'_k=\mSi$.
\end{proof}

The most interesting implication is probably $(ii) \Leftrightarrow (v)$. In particular, it means that removing from $S$ or adding to $S$ an active edge does not change the types of edges, nor \emph{a fortiori} the set of active edges. This can be also be rewritten in term of subgraph intervals.

\begin{cor} \label{cor:inter}
For any subgraph $S$, the equivalence class of $S$ is exactly the subgraph interval \mbox{$[S \backslash \Act(S), S \cup \Act(S)]$}.
\end{cor}




\subsection{Indexation of the intervals by spanning trees}

\begin{theo} \label{partition}
Let $G$ be a graph and $\Delta$ a decision tree. Then the set of subgraphs can be partitioned into subgraph intervals indexed by the spanning trees:
\begin{equation} 2^{E(G)} = \biguplus_{T\textrm{ spanning tree of }G} [T \backslash \mathcal I(T), T \cup \mathcal E(T)], \label{eqpart} \end{equation}
where $\mathcal I(T)$ and $\mathcal E(T)$ are the sets of internal and external $\Delta$-active edges of the spanning tree $T$.
\end{theo}

Note  that Corollary \ref{cor:inter} tells that each of the intervals $[T \backslash \mathcal I(T), T \cup \mathcal E(T)]$ constitutes an equivalence class for the relation $\sim$.

We  first prove that every subgraph of $G$ is equivalent to some spanning tree.

\begin{lem}
\label{reptree}
For every $S$ subgraph of $G$, the set of edges with type \bSi \,or \bI\, forms a spanning tree of $G$ equivalent to $S$.
\end{lem}
\begin{proof} Let us denote by $T$ the set of edges with type \bSi\, and \bI. We choose the version of Algorithm \ref{type} where edges of type \bL\, are deleted and edges of type \bI\, are contracted. Let us prove that $T$ has no cycle and is connected, which exactly means that $T$ is a spanning tree.

\textbf{1. $\boldsymbol T$ has no cycle.} Let us assume that $T$ has a cycle $C$. Consider $e_k$ the maximal edge of $C$ for the $(\Delta,S)$-ordering.   At the $k$-th iteration of the algorithm, every edge of $C$ is contracted except $e_k$, because $C$ is only made of edges of type \bSi\, or \bI. This implies that $e_k$ is a loop, which contradicts the fact that $e_k$ has type \bSi\, or \bI.

\textbf{2. $\boldsymbol T$ has only one connected component} (This is precisely the dual of the previous point.). Let us assume that $T$ has more than one connected component. This means that there exists a cocycle $D$ of $G$ only made of edges of type \bSe\, and \bL. Consider $e_k$ the maximal edge of $D$ for the $(\Delta,S)$-ordering. Then at the $k$-th iteration of the algorithm, every edge of $D$ has been deleted except $e_k$. This implies that $e_k$ is a isthmus, which contradicts the fact that $e_k$ has type \bSe\, or \bL.

\textbf{3. $\boldsymbol T$ is equivalent to $\boldsymbol S$.} We have $\diffs S T \subseteq \Act(S)$ because the edges of type \bSi\, (resp. \bSe) for $S$ are internal (resp. external) in both subgraphs. By the implication (iv)~$\Rightarrow$~(i) of Proposition \ref{samehistory}, this means that $S \sim T$. 
\end{proof}

\begin{proof}[Proof of Theorem \ref{partition}]
Since the equivalence classes for the relation $\sim$ partition $2^{E(G)}$, we only need to show that there exists a unique spanning tree inside each of these classes. The existence is proved by Lemma \ref{reptree}. Now let us show the uniqueness of the spanning tree.

Let $T$ and $T'$ be two spanning trees of $G$ such that $T \sim T'$. Thus they share the same partition of $E(G)$. But remember that edges of type \bSi\, are always internal and those of type \bSe\, always external. Furthermore,  Proposition \ref{intextact} tells that each edge of type \bI\, is internal and each edge of type \bL\, is external. Therefore $T = T'$.
\end{proof}
\subsection{Some descriptions of the Tutte polynomial}


Let us fix a graph $G$ and a decision tree $\Delta$. 
The following lemma describes how the number of connected components and the cyclomatic number behave when we add/remove an active edge.

\begin{lem} 
\label{addremove}
Let $S$ be a subgraph of $G$ and $e$ an active edge for $S$.
\begin{enumerate}
\item[(a)] If $e$ is external and has type \bL, then
$$\cc(S \cup \{ e \})= \cc(S),  \quad
\cycl(S \cup \{ e \})= \cycl(S) + 1.$$
\item[(b)] If $e$ is internal and has type \bI, then
$$\cc(S \backslash  e )= \cc(S) + 1,  \quad
\cycl(S \backslash  e )= \cycl(S). $$
\end{enumerate}
\end{lem}

\begin{proof}
\textbf{(a) $\boldsymbol e$ is external and has type \bL.} Point (1) from Lemma \ref{cyclecocycle} ensures that there exists a path only made of edges of type \bSi\, (so this is a path in $S$) linking the endpoints of $e$. Therefore including $e$ in $S$ does not add a connected component, hence $\cc(S \cup \{ e \})= \cc(S)$. Moreover,
$$\cycl(S \cup \{ e \})= |S \cup \{ e \}| + \cc(S \cup \{ e \}) - |V(G)| =   |S| + 1 + \cc(S) - |V(G)| =  \cycl(S) + 1.$$

\textbf{(b) $\boldsymbol e$ is internal and has type \bI.} Point (2) from Lemma \ref{cyclecocycle} ensures that there exists a cocycle only made of $e$ and edges of type \bSe. In other terms, there exists no path with edges of $S \backslash e$ linking the endpoints of $e$. So removing $e$ from $S$ will increase the number of connected components: $\cc(S \backslash  e )= \cc(S) + 1.$ Moreover,
$$\cycl(S \backslash  e )= |S \backslash  e | + \cc(S \backslash  e ) - |V(G)|  =   |S| - 1 + \cc(S) + 1 - |V(G)| =  \cycl(S),$$
which ends the proof. \end{proof}

\begin{prop} \label{rouge}
Let $T$ be a spanning tree. We denote by $I_T$ the set of the subgraphs equivalent to $T$. Then
\begin{equation} 
\label{eq:weight}
\sum_{S \in I_T}(x-1)^{\cc(S)-1}(y-1)^{\cycl(S)} = x^{|\mathcal I(T)|} y^{|\mathcal E(T)|},
\end{equation}
where $\mathcal I(T)$ (resp. $\mathcal E(T)$) denotes the set of internal (resp. external) $\Delta$-active edges of $T$.
\end{prop}

\begin{proof}
By Corollary \ref{cor:inter} the interval $I_T$ is the set of subgraphs of the form  $(T \backslash R_i) \cup R_e$, where $R_i \subseteq \mathcal I(T)$ and $R_e \subseteq \mathcal E(T)$.
In other terms, each subgraph of $I_T$ can be obtained by deleting from $T$ some edges of type \bI\, one by one, then adding to $T$ some edges of type \bL. Thus, by repeated applications of Lemma \ref{addremove}, given $S=(T \backslash R_i) \cup R_e \in I_T$, the number $\cc(S)$ is equal to $\cc(T) + |R_i| = |R_i| + 1$ and $\cycl(S)$ is equal to $\cycl(T) + |R_e| = |R_e|$. The formula \eqref{eq:weight} then derives from the identity $(X+1)^{|A|}(Y+1)^{|B|}=\sum_{\substack{S_A \subseteq A \\ S_B \subseteq B}} X^{|S_A|}Y^ {|S_B|}.$
%
%
%
%
\end{proof}

\noindent \textbf{Remark.} Theorem \ref{charact} can be easily deduced from the previous proposition. Indeed, the formula \eqref{celuila} can be obtained by summing \eqref{eq:weight} over all spanning trees $T$. (Remember that the intervals $I_T$ form a partition of the subgraph set -- see Theorem \ref{partition}.)

We can adapt the previous reasoning to obtain more descriptions of the Tutte polynomial.

\begin{prop} Given any decision tree $\Delta$, the Tutte polynomial of $G$ admits the three following descriptions:
\begin{equation}
\label{eqfor}
T_G(x,y) = \sum_{F\textrm{ spanning forest of $G$}} (x-1)^{\cc(F)-1}y^{\ell(F)},
\end{equation}
\begin{equation}
\label{eqfive}
T_G(x,y) = \sum_{K\textrm{ connected subgraph of $G$}} x^{i(K)}(y-1)^{\cycl(K)},
\end{equation}
\begin{equation}
\label{eqsix}
T_G(x,y) = \sum_{S\textrm{ spanning subgraph of $G$}} \left(\frac x 2 \right)^{i(S)} \left(\frac y 2 \right)^{ \ell(S) },
\end{equation}
where $i(S)$ (resp. $\ell(S)$) denotes the number of edges with type \bI\, (resp. \bL) for a subgraph $S$.
\end{prop}

These three descriptions can be also found in \cite{gordon-traldi}, but in the specific case of ordering activity.

\begin{proof} \textbf{Formula  \eqref{eqfor}.} Let $T$ be a spanning tree and $I_T$ the set of the subgraphs equivalent to $T$. Looking back on the proof of Proposition \ref{rouge}, we see that the spanning forests of $I_T$ are of the form $T \backslash R_i$ with $R_i \subseteq \mathcal I(T)$. So every subgraph of $I_T$ can be written as $F \cup R_e$, where $F$ is a spanning forest of $I_T$ and $R_e \subseteq \mathcal E(T)$. By Proposition \ref{samehistory}, the set $\mathcal E(T)$, which consists of edges with type \bL\, for $T$, equals $\mathcal L(F)$ for any forest $F$ in $I_T$, where $\mathcal L(F)$ denotes the set of edges with type \bL\, for $F$. Hence, we have 
$$I_T = \biguplus_{F\textrm{ spanning forest of }I_T} \left[F,F \cup \mathcal L(F)\right].$$
As the intervals $I_T$ partition the set of subgraphs (see Theorem \ref{partition}), we deduce:
\begin{equation} \label{equite}
2^{E(G)} = \biguplus_{F\textrm{ spanning forest of }G} [F, F \cup \mathcal L(F)].
\end{equation}
A subgraph $S \in [F, F \cup \mathcal L(F)]$ is of the form  $F \cup R$, with $R \subseteq \mathcal L(F)$, and thanks to Lemma \ref{addremove} satisfies $\cc(S) = \cc(F)$ and $\cycl(S) = \cycl(F) + |R| = |R|$.
So we get
$$\sum_{S \in [F, F \cup \mathcal L(F)]} (x-1)^{\cc(S)-1}(y-1)^{\cycl(S)} = (x-1)^{\cc(F)-1} y^{|\mathcal L(F)|}.$$
The formula \eqref{eqfor} results from the previous equation and \eqref{equite}.

\noindent \textbf{Formula  \eqref{eqfive}.} Similar to the previous point, in a dual way.
 
\noindent \textbf{Formula  \eqref{eqsix}.} Let $T$ be a spanning tree. By Proposition \ref{samehistory}, each subgraph $S$ equivalent to $T$ satisfies $i(S)=i(T)$ and $\ell(S)=\ell(T)$. So the sum   $\sum_{S \sim T} \left(\frac x 2 \right)^{i(S)} \left(\frac y 2 \right)^{ \ell(S) }$ equals $\left(\frac x 2 \right)^{i(T)} \left(\frac y 2 \right)^{ \ell(T)}$ times the number of subgraphs equivalent to $T$, namely $2^{|\Act(T)|}=2^{i(T)+\ell(T)}$. So the previous sum is equal to $x^{i(T)}y^{\ell(T)}$. Summing over all spanning tree $T$, we obtain formula \eqref{eqsix}. \end{proof}

\section{Back to the earlier activities}
\label{sec:spec}

It is time to show that each activity defined at Section \ref{sec:activity} is  a $\Delta$-activity. In this way, we prove that all these activities are Tutte-descriptive.
The key tools for this section are Theorem \ref{treecompatible} and Corollary \ref{cor:act}. 

\subsection{Ordering activity}

Let $G$ be a graph and consider a linear order $<_{ord}$:
$\epsilon_1 < \epsilon_2 < \dots < \epsilon_m,$
where $E(G) = \ens{\epsilon_1,\dots,\epsilon_m}$.
Define the decision function $\Delta_{ord}$ by setting
\begin{equation} \Delta_{ord}(d_1,\dots,d_k) =  \epsilon_{m-k} \end{equation}
for every sequence of directions $(d_1,\dots,d_k)$ with $k < m$.

\begin{prop}
The ordering activity is equal to the  $\Delta_\phi$-activity. Therefore it is Tutte-descriptive.
\end{prop}

We recall that the ordering activities (also called Tutte's activities) have been introduced in Subsection~\ref{ss:ord}.

\begin{proof} Fix a spanning tree $T$. It is straightforward to see from the definition of $\Delta_{ord}$ that $\epsilon_m < \dots < \epsilon_1$
for the $(\Delta_{ord},T)$-order. In other terms, the $(\Delta_{ord},T)$-order is inverse to $<_{ord}$. So by Proposition~\ref{maximal}, an edge is ordering-active if and only if it is $\Delta_{ord}$-active. We conclude thanks to Theorem~\ref{charact}, applied to the decision function $\Delta_{ord}$.
\end{proof}

\noindent \textbf{Remark.} Alternatively, we could have also used Corollary \ref{cor:act} with constant order map $T \mapsto\  >_{ord}$.

\begin{figure}[h!]
\begin{center}
\includegraphics[scale=1]{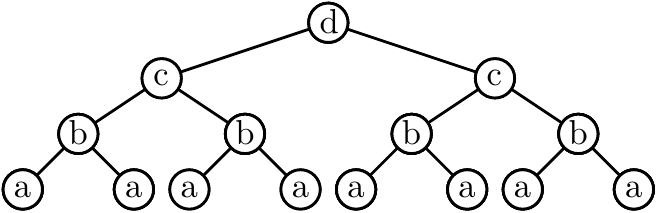}
\end{center}
\caption{The decision tree corresponding to the linear order $a < b < c < d$.}
\label{fig:dttut}
\end{figure}

\noindent  \textbf{Example}. The decision tree associated to the graph from Figure \ref{fig:extut} with the linear order $ a < b < c < d$  is depicted at Figure \ref{fig:dttut}.

\subsection{Embedding activity}

Let $G$ be a graph. We want to express the embedding activities (see Subsection \ref{ss:bernardi} for the definition) as $\Delta$-activities and thus give an alternative proof of Theorem 7 from \cite{bernardi-tutte}.

\begin{prop}
\label{prop:ber}
For any embedding $M_G$ of the graph $G$, the embedding activity is a  $\Delta$-activity and so a Tutte-descriptive activity. 
\end{prop}

\begin{figure}[h!]
\begin{center}\includegraphics[scale=1]{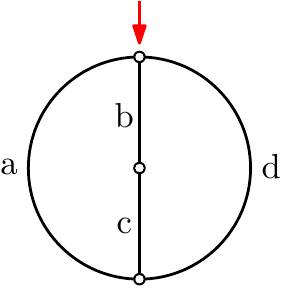} \hspace{2cm}
\includegraphics[scale=1.2]{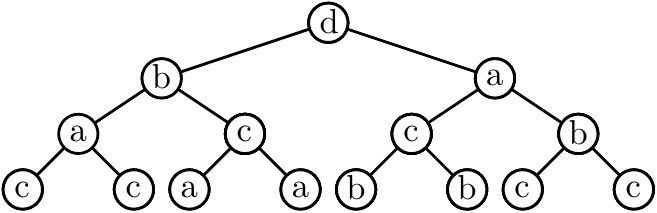}
\end{center}
\caption{A map and a decision tree for the embedding activity.}
\label{fig:dqber}
\end{figure}

\noindent \textbf{Example. }Let $M$ and $\Delta$  be the map and the decision tree of  Figure \ref{fig:dqber}. Consider $T = \ens{b,d}$. The $(M,T)$-ordering \footnote{Reminder: the $(M,T)$-order is the order of first visit during the tour of $T$ -- see Subsection \ref{ss:bernardi}.} is $a < b < c < d$. So the set of $(M,T)$-active edge is $\ens{a,b}$. The $(\Delta,T)$-ordering is $d < a < c < b$, so the set of $\Delta$-active edges for $T$ is as well $\ens{a,b}$. We can check that the two sets of active edges coincide (even if the two orderings do not look related).

 The proof is in two steps, embodied in two lemmas. The first lemma uses tree-compatibility via Corollary~\ref{cor:act}.

\begin{lem}
For each embedding $M_G$ of a graph $G$, the activity that maps any spanning tree $T$ onto the set of edges that are \emph{maximal} for the ($M_G$,$T$)-order in their fundamental cycles/cocycles is a $\Delta$-activity and so Tutte-descriptive.
\end{lem}

\begin{proof} It suffices to prove that the order map that sends a spanning tree $T$ of $G$ onto the $(M_G,T)$-order satisfies the hypothesis of Corollary \ref{cor:act}.

For any spanning tree $T$, we denote by $\phi_i(T)$ (resp. $h_i(T)$) the $i$-th smallest edge (resp. half-edge) for the $(M_G,T)$-ordering, namely the edge (resp. half-edge) which is visited in $i$-th position during the tour of the tree $T$.  We consider  $T$ and $T'$ two spanning trees of $G$ and $k$ an integer from $\ens{0,\dots,|E(G)|-1}$ such that \begin{equation} T \cap \ens{\phi_1(T),\dots,\phi_k(T)} = T' \cap \ens{\phi_1(T),\dots,\phi_k(T)}.
\label{bertc}
\end{equation} Moreover, let $\ell$ be the integer in $\ens{1,\dots,2|E(G)|}$ such that $h_\ell(T)$ is the smallest half-edge of $\phi_{k+1}(T)$. We want to show that $\phi_i(T) = \phi_i(T')$  for each $i \in \ens{1, \dots, k+1}$.

Let us prove by recurrence on $j \in \ens{1,\dots,\ell}$ that $h_j(T) = h_j(T')$. For $j=1$, the half-edges $h_1(T)$ and $h_1(T')$ are both equal to the root half-edge of $M_G$. Now assume that for some $j \in \ens{1,\dots,\ell-1}$ we have $h_j(T) = h_j(T')$. Let $\phi_i(T)$ be the edge that corresponds to $h_j(T)$. The edge $\phi_i(T)$ belongs to $\ens{\phi_1(T),\dots,\phi_k(T)}$  since we have $h_j(T) < h_{\ell}(T)$ for the $(M_G,T)$-ordering. So the equivalence 
$$\phi_i(T)\textrm{ internal for }T \ \Leftrightarrow \phi_i(T)\textrm{ internal for }T',$$
that results from \eqref{bertc},
holds. Thus, we have $$t(h_j(T);T) = t(h_j(T');T')$$
(let us recall that $t(.;T)$ denotes the motion function of the spanning tree $T$, see Subsection \ref{ss:bernardi}), that is 
$$h_{j+1}(T) = h_{j+1}(T'),$$
which proves the recurrence.

We have shown that whether we consider the tour of $T$ or the tour of $T'$, the first visited half-edges are $h_1(T),\dots,h_\ell(T)$ in this order. Consequently, the $k+1$ first visited edges are identical in both cases, namely $\phi_1(T),\dots,\phi_k(T),\phi_{k+1}(T)$ in this order. This means that for every $i \in \ens{1,\dots,k+1}$ we have $\phi_i(T) = \phi_i(T')$.
\end{proof}

Despite the similarities with embedding activities, the previous lemma is not sufficient to recover Bernardi's activity, for which an active edge is \emph{minimal} in its fundamental cycle/cocycle rather than \emph{maximal}. 

We could then have the idea to adapt the previous proof by considering the order map that sends each spanning tree $T$ onto the \emph{reversed} ($M_G$,$T$)-ordering. But a problem occurs: this order map is not tree-compatible\footnote{Indeed, when an order map $\phi$ is tree-compatible, the minimal edge for $\phi(T)$ is the same for all spanning trees $T$, which is obviously not the case here.}.

\begin{figure}[h!]
\begin{center}
\includegraphics[scale=1]{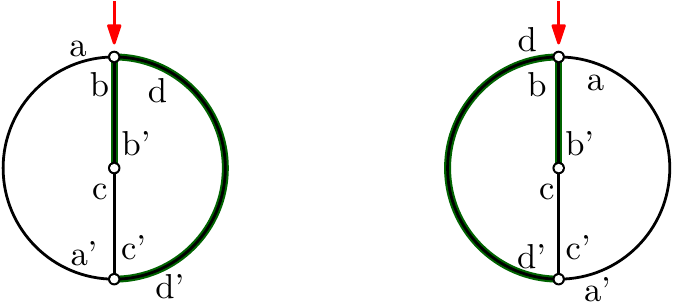}
\end{center}

\caption{A map $M$ with a spanning tree $T$ and its mirror map $M'$.}
\label{fig:mirror}
\end{figure}

We could also have the seemingly desperate idea to "reverse" the map instead of the ($M_G$,$T$)-ordering. For example, consider the rooted maps $M$ and $M'$ equipped with the spanning tree $T$ of Figure \ref{fig:mirror}. The map $M'$ is the mirror map of $M$. The $(M,T)$-ordering for half-edges equals $a < b < c < b' < d < c' < a' < d'$ while the $(M',T)$-ordering equals $a > b' > c > b > d' > c' > a' > d$\footnote{If we formally remove the primes in these two orderings, we will observe a strange phenomenon: the two orderings are reverse! An explanation is implicitly given in the proof of Lemma \ref{miror}.}. The orderings on the edges are not reverse ($\ar a < \ar b < \ar c < \ar d$ for $M$ and $ \ar d < \ar a < \ar c < \ar b$ for $M'$) but the $(M,T)$-active edges correspond to the edges that are maximal in their fundamental cycle/cocycle for the $(M',T)$-ordering, namely $\ar a$ and $\ar b$. (For instance, take $\ar a$: its fundamental cycle is $\ens{\ar a,\ar d}$. We have $\ar a < \ar d$ for the $(M,T)$-order and $\ar d < \ar a$ for the $(M',T)$-order.) It turns out that this property is general. 

\begin{lem} \label{miror}
Consider any embedding $M_G= (H,\sigma,\alpha)$ of the graph $G$ rooted on a half-edge $a$, and denote by $M^\#_G$ the mirror map of $M_G$, that is to say the map $(H,\sigma^{-1},\alpha)$ rooted on the half-edge $\sigma^{-1}(a)$.

For any spanning tree $T$, an internal (resp. external) edge is $(M_G,T)$-active if and only if it is maximal in its fundamental cocycle (resp. cycle) for the $(M^\#_G,T)$-ordering. 
\end{lem}

This lemma shows that if we change "minimal" by "maximal" in the definition of the embedding activity, we obtain a strictly equivalent notion (and maybe more natural).

\begin{proof} Fix $T$ a spanning tree. We denote by $t$ the motion function corresponding to $M_G$ and $T$, and by $t'$ the motion function corresponding to $M^\#_G$ and $T$.

\textbf{1.} Let us prove that an edge $\ar h$ is smaller than $\ar g$ for the $(M^\#_G,T)$-ordering if and only if we have $\max(g,g') < \max(h,h')$ for the  $(M_G,T)$-ordering. Let $\ar g$ be an edge and $m$ the number of edges. Since $t'$ is a cyclic permutation of the set of half-edges (see Lemma \ref{tourcyclique}), there exist $k$ and $\ell$ in $\ens{0,\dots,2m-1}$ such $g = t'^k(\sigma^{-1}(a))$ and $g' = t'^\ell(\sigma^{-1}(a))$.
Moreover, it is easy to see from the algebraic definition of the  motion function that $t' = \sigma^{-1} \circ t^{-1} \circ \sigma$ and that $\sigma^{-1}(t(\hat g)) \in \ens{\hat g,\alpha(\hat g)}$ for every half-edge $\hat g$. Hence, 
$$g=t'^k(\sigma^{-1}(a)) = \sigma^{-1}(t^{-k}(a)) = \sigma^{-1}(t^{2m-k}(a)) \in \ens{t^{2m-k-1}(a),\alpha(t^{2m-k-1}(a))}.$$
Similarly, we have $g' \in \ens{t^{2m-\ell-1}(a),\alpha(t^{2m-\ell-1}(a))}$. But $g$ and $g'$ belong to the same edge, so 
$$\ar g =  \ens{t^{2m-k-1}(a), \alpha(t^{2m-k-1}(a))} = \ens{t^{2m-\ell-1}(a), \alpha(t^{2m-\ell-1} (a))}$$
and since $k \neq \ell$ (we have $g \neq g'$), we deduce that 
$$ \ar g =  \ens{t^{2m-k-1}(a),t^{2m-\ell-1}(a)}.$$
Thus, the following two equalities hold:
$$\min(g,g')= \min \left(t'^k(\sigma^{-1}(a)),t'^\ell(\sigma^{-1}(a))\right) = t'^{\min(k,\ell)}(\sigma^{-1}(a)),$$  where the $\min$  concern the $(M^\#_G,T)$-ordering,
and
$$ \max(g,g')= \max \left(t^{2m-k-1}(a),t^{2m-\ell-1}(a)\right) = t^{2m - \min(k,\ell) - 1}(a), $$ where the $\max$  concern the $(M^\#_G,T)$-ordering,
Let $\ar h$ be an edge different from $\ar g$ with $h = t'^i(\sigma^{-1}(a))$ and $h' = t'^j(\sigma^{-1}(a))$. Using the above equalities, there is equivalence between the following statements:
$$\begin{array}{lcl}
\ & \ar h < \ar g &  \textrm{ for the }(M^\#_G,T)\textrm{-ordering}, \\
\Leftrightarrow \  & \min (h,h') < \min (g,g') &  \textrm{ for the }(M^\#_G,T)\textrm{-ordering}, \\
\Leftrightarrow \  & \min(i,j) < \min(k,\ell) &  \\
\Leftrightarrow \ & 2m - \min(k,\ell) - 1 < 2m - \min(i,j) - 1 &  \\
\Leftrightarrow \ & \max(g,g') < \max(h,h') &  \textrm{ for the }(M_G,T)\textrm{-ordering}.
\end{array}$$

\textbf{2.} Let $\ar g$ be an external (resp. internal) edge and $\ar h$ an edge in the fundamental cycle (resp. cocycle) of $\ar g$ with $g < g'$, $h < h'$ and $g < h$, where $<$ denotes from now on the $(M_G,T)$-ordering. Let us show that $g' < h'$.

Using Lemma 4 from \cite{bernardi-tutte}, it is not hard to see that for every spanning tree $T$, deleting an external edge or contracting an internal edge of $M_G$ does not change the $(M_G,T)$-ordering between the remaining half-edges.
This is why we can restrict the proof of this point to the case where $E(G) = \ens{\ar g,\ar h}$. As it must contain a cycle  with 2 edges or a cocycle with 2 edges, $G$ must consist of two vertices linked by $\ar g$ and $\ar h$. Since $\ar g < \ar h$, the root of the map is $g$. The only two possibilities for $(G,T)$ are depicted in Figure \ref{baseber}.
Thus, we must have in both cases $t = (g,h,g',h')$.  Therefore $g' < h'$.

\fig{[scale=1]{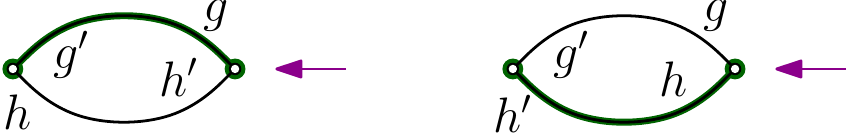}}{Two possibilities when $G$ is restricted to two edges}{baseber}

\textbf{3.} Fix   a $(M_G,T)$-active edge $\ar g$. By definition, $\ar g$ is minimal for the $(M_G,T)$-ordering in its fundamental cycle/cocycle $C$. Let $\ar h$ be any other edge in $C$. If we assume $g < g'$ and $h < h'$, then $g < h$. By point 2, we have $g' < h'$. So by point 1, $\ar g$ is greater than $\ar h$ for the $(M^\#_G,T)$-ordering. This being true for all $\ar h$, the edge $\ar g$ is maximal for the $(M^\#_G,T)$-ordering inside $C$.

\textbf{4.} Conversely, fix an edge $\ar g$ that is maximal for the $(M^\#_G,T)$-ordering in its fundamental cycle/cocycle $C$, with $g < g'$. For any other edge $\ar h$ in $C$ with $h < h'$, we have $g' < h'$ by point $1$. By contraposition of point 2, we have $g < h$. This means that inside its fundamental cycle/cocycle, $\ar g$ is minimal for the $(M_G,T)$-ordering, and thus $(M_G,T)$-active.\end{proof}

Finally, the combination of the two previous lemmas gives a proof of Proposition~\ref{prop:ber}, since $M^\#_G$ is an embedding of the graph $G$.

\subsection{Blossoming activity}

Now we are going to deal with the blossoming activity of Subsection \ref{ss:blo}. We fix an embedding $M$ of $G$.
We begin by a lemma (with a $\Delta$-activity flavour) establishing another characterization of internal blossoming-active edges.


\begin{lem}
For any spanning tree $T$ of $G$, an internal edge $e$ is blossoming-active if and only if $e$ is an isthmus of $M'$ when it is visited for the first time during the computation of $\tau(T)$ (cf. Algorithm \ref{prune}).
\end{lem}
\begin{proof} 
%
Let us compare the executions of Algorithm \ref{prune} with inputs $T$ and $T \backslash e$. Before the first visit of $e$, the map $M'$ is the same. At this point, there are two possibilities for $e$.

If $e$ is not an isthmus in $M'$, then $e$ will be deleted when the input is $T \backslash e$. But when the input is $T$, the edge $e$ will be never deleted since it is internal. So $e \notin \tau(T \backslash  e)$ and $e \in \tau(T)$. Thus $\tau(T \backslash e) \neq \tau(T)$: the edge $e$ is not blossoming-active.

If $e$ is an isthmus in $M'$, then $e$ will not be deleted in both cases, and this for the rest of the run. Each other edge having the same status internal/external for the two inputs, we will have $\tau(T \backslash e) = \tau(T)$: the edge $e$ is blossoming-active.
\end{proof}

We can now prove that our internal edge activity can be extended into a Tutte-descriptive activity.

\begin{prop} \label{prop:blo}
For any embedding of the graph $G$, the internal blossoming activity can be extended into a $\Delta$-activity and so a Tutte-descriptive activity.
\end{prop}

\begin{figure}[h!]
\begin{center}
\includegraphics[scale=1]{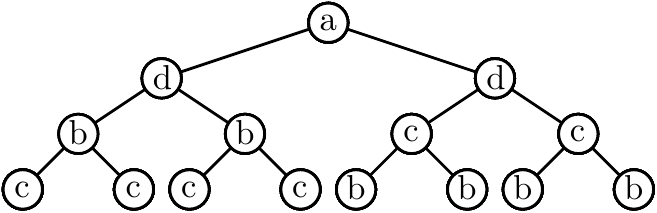}
\end{center}
\caption{A decision tree for the blossoming activity, for  the embedded graph from Figure \ref{fig:exblo}.}
\label{fig:dtblo}
\end{figure}

\textbf{Example.} Consider the map from Figure \ref{fig:exblo}. A suitable decision tree $\Delta$ is shown in Figure \ref{fig:dtblo}. For instance, consider the spanning tree $T = \ens{c,d}$. The $(\Delta,T)$-ordering is $a < d < b < c$ and the only internal active edge is $c$, as for the blossoming activity.

\begin{proof}
For any spanning tree $T$, let $\phi(T)$ denote the first visit order of the edges of $G$ in a run of Algorithm \ref{prune} with input $T$. For example, if we consider the first embedded graph from Figure \ref{fig:exblo}, then we map the spanning tree $\ens{c,d}$ onto the ordering $a < d < b < c$. Moreover, let $\phi_k(T)$ be the $k$-th smallest edge of $G$ for $\phi(T)$. We want to use Theorem~ \ref{treecompatible}.

Consider $T$ and $T'$ two spanning trees of $G$ and $k \in \ens{0,\dots,|E(G)|-1}$ such that $$T \cap \ens{\phi_1(T),\dots,\phi_k(T)} = T' \cap \ens{\phi_1(T),\dots,\phi_k(T)}.$$
In Algorithm \ref{prune}, only the status (external, internal, isthmus) of $e$ in $M'$ at each iteration has an influence on the next values of $e$, $h$ and $M'$. But before the visit of $\phi_{k+1}(T)$, the statuses of $e$ are the same in $T$ and  in $T'$, since $T \cap \ens{\phi_1(T),\dots,\phi_k(T)} =$ \mbox{$T' \cap \ens{\phi_1(T),\dots,\phi_k(T)}$}.
So we have $\phi_i(T) = \phi_i(T')$ for every $i \in \ens{1,\dots,k+1}$.  

By Theorem \ref{treecompatible}, the order map $\phi$ is tree-compatible, meaning that we can construct a decision tree $\Delta_\phi$ such that the $(\Delta_\phi,T)$-ordering coincides with $\phi(T)$. In other terms, 
the edges are visited in the same order in Algorithm \ref{prune} and in Algorithm \ref{type} (with decision tree $\Delta_\phi$).

By Theorem \ref{charact}, the activity that sends a spanning tree onto the set of its $\Delta_\phi$-active edges is Tutte-descriptive. Let us show that it extends the internal blossoming activity. Algorithm~\ref{algint}, when executed on a spanning tree $T$ with decision tree $\Delta_\phi$, outputs the set of edges that are isthmuses at the time of their first visit in Algorithm~\ref{prune}. By Proposition \ref{prop:algint}, this is the set of internal $\Delta_\phi$-active edges. But by the previous lemma, this coincides also with the set of internal blossoming-active edges.
\end{proof}

This proof allows us to give a natural definition of the complete blossoming activity: an external edge for a spanning tree $T$ is  \textit{blossoming-active} if it is $\Delta_\phi$-active, where $\Delta_\phi$ is any decision tree compatible with the order map $\phi$ defined in the previous proof. The $\Delta_\phi$-active edges are uniquely determined, although $\Delta_\phi$ is not. Indeed, we can give an intrinsic characterization of the (complete) blossoming activity, like the following one.

\begin{cor} 
\label{corblo}
Given any spanning tree $T$, an edge is blossoming-active if and only if it has been visited last in its fundamental cycle/cocycle during the run of Algorithm \ref{prune} with input $T$.
\end{cor}
\begin{proof} By definition of $\Delta_\phi$, the $(\Delta_\phi,T)$-ordering corresponds to the order of first visit during the run of Algorithm \ref{prune} with input $T$. We conclude using Proposition \ref{maximal}.
\end{proof}

Finally, let us give a description of the preimage of a spanning tree $T$ under $\tau$. (We recover a property similar to what we have already seen for the DFS activity.)

\begin{prop} \label{prop:tauf}
For each spanning tree $T$ and each spanning forest $F$, we have
$\tau(F) = T$ if and only if $F \in [T \backslash \mathcal I (T), T]$  (see Subsection \ref{sss:interval} for the definition of an interval),
where $\mathcal I$ denotes the internal blossoming activity.
\end{prop}
\begin{proof} Given the spanning forest $F$, the edges we delete in Algorithm \ref{prune} are the external edges that are not isthmuses at the time of their first visit. They are precisely the edges of $\Delta_\phi$-type \bSe\, or \bL\, for $F$. Thus, $\tau(F)$ corresponds to the set of edges of type \bSi\, or \bI\, for $F$. So by Lemma \ref{reptree}, the spanning tree $\tau(F)$ is equivalent to $F$. By Corollary \ref{cor:inter}, we have then $F \in [\tau(F) \backslash \mathcal I (\tau(F)), \tau(F) \cup \mathcal E(\tau(F))]$, where $\mathcal E$ is the external blossoming activity. 
So, by Theorem \ref{partition}, we have $F \in [T \backslash \mathcal I (T), T \cup \mathcal E(T)]$ if and only if $T = \tau(F)$. But the restriction of the interval $[T \backslash \mathcal I (T), T \cup \mathcal E(T)]$ to the spanning forests of $G$ is $[T \backslash \mathcal I (T), T]$: indeed, Lemma \ref{addremove} states that a subgraph with an internal edge of type \bL\, has at least a cycle.
\end{proof}

\subsection{DFS activity}
\label{ss:dfst}

We end this section with DFS activity defined in Subsection \ref{ss:dfs}. Let us recall that we now consider a graph $G$ without multiple edges.

\begin{prop}
\label{dfstut}
For any labelling of $V(G)$ with integers $1,2,\dots,|V(G)|$, the external DFS-activity restricted to spanning trees can be extended into a $\Delta$-activity and so a Tutte-descriptive activity.
\end{prop}

\begin{figure}[h!]
\begin{center}
\includegraphics[scale=1]{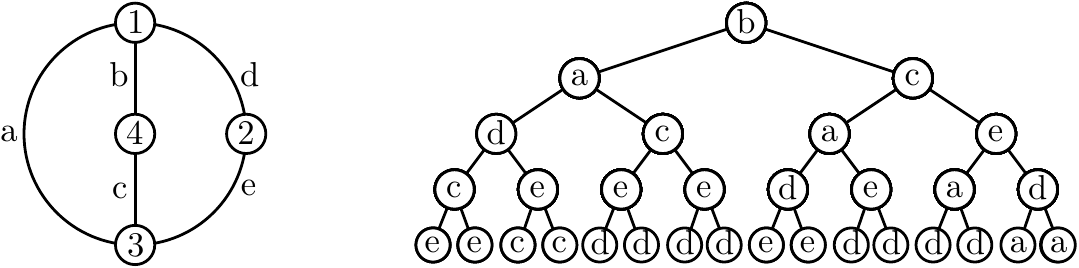}
\end{center}
\caption{A graph with no multiple edges and the decision tree corresponding to the DFS-activity.}
\label{fig:dtdfs}
\end{figure}

\noindent \textbf{Example. } We show in Figure \ref{fig:dtdfs} a graph with a decision tree $\Delta$ inducing the DFS activity. For instance, consider the spanning tree $T = \ens{a,c,e}$. The $(\Delta,T)$-ordering is $b < a < c < e < d$, so the only external $\Delta$-active edge is $d$. One can check that $d$ is also the only external DFS-active edge.

\begin{proof} The principle of the proof can be detailed as follows:
\begin{itemize}
\item we define an order map $\phi$ that is related to the visit order in the (greatest-neighbor) DFS, which is described by Algorithm \ref{DFS}
\item we prove that $\phi$ is tree-compatible by checking the hypotheses of Corollary \ref{cor:act},
\item we show that the external DFS-activity can be described in terms of maximality in $\phi$.
\end{itemize}

We would want to match an order map $\phi$ with the visit order in Algorithm \ref{DFS}. But this algorithm only considers the internal edges while we also need to order the external edges to define an order map. That is why we need to enrich Algorithm \ref{DFS} to take into account external edges. The result is Algorithm \ref{pmDFS}. 
\begin{algorithm}[h!]
\caption{The order map for DFS activity }

\label{pmDFS}
\begin{algorithmic}[5]
\Require {\color{darkgray}$H$ spanning subgraph of $G$}.
\Ensure {\color{darkgray}A total ordering of $E(G)$ given by the edges $\phi_1(H),\phi_2(H),\dots,\phi_{|E(G)|}(H)$ in this order.}
\State $\mathcal F(H) \leftarrow \emptyset$;
\State $j \leftarrow 0$; 
\While {there is a unvisited vertex} 
	\State  \textbf{mark} the least unvisited vertex of $G$ \textbf{as visited};
	\While {{\color{darkgray}there is a visited vertex with unvisited incident edges}}
		\State \ \Comment{In this algorithm "incident" means "incident in $G$" (not in $H$).}
		\State $v \leftarrow$ the most recently visited such vertex;
		\While{{\color{darkgray} $v$ has an unvisited incident edge}}
			\State $u \leftarrow$ the greatest neighbor of $v$ {\color{darkgray}linked in $G$ by an unvisited edge};
			\State $j \leftarrow j + 1$;
			\State $\phi_j(H) \leftarrow \ens{u,v}$;
			\State {\color{darkgray}\textbf{mark} $\phi_j(H)$ \textbf{as visited}};
			\If{{\color{darkgray}$\phi_j(H)$ is internal \textbf{and} if $u$ is unvisited}} \label{lineuse}
				\State \textbf{mark} $u$ \textbf{as visited};
				\State $v \leftarrow u$; 
				\State \textbf{add} $\phi_j(H)$ in $\mathcal F(H)$;
			\EndIf
		\EndWhile
	\EndWhile
\EndWhile 
\State \Return {\color{darkgray}$\phi_1(H) < \phi_2(H) < \dots < \phi_{|E(G)|}(H)$};
\end{algorithmic}
\end{algorithm}

\noindent \textbf{Informal description.} The input is a subgraph $H$ of $G$. As in Algorithm \ref{DFS}, we begin by the least vertex. We proceed to the DFS of the graph $\boldsymbol G$ that favors the largest neighbors with an extra rule: when we visit an external edge (for $H$), we do not go through it, we come back to the original vertex as if this edge had never existed. Moreover, the visited edges (internal and external) are marked so that we visit them only once each. The rest of the algorithm is exactly as Algorithm~\ref{DFS}. The output is the visit order of the edges, denoted  by $\phi(H)$, instead of the DFS forest. The changes between Algorithm \ref{DFS} and \ref{pmDFS} are indicated in gray.

\noindent \textbf{Remark.} In this proof, only connected subgraphs are important. We could have simplified Algorithm \ref{pmDFS} with a  restriction of the input to connected subgraphs, but it could have interfered with understanding. 

\noindent \textbf{Example.} A run of Algorithm \ref{pmDFS} is depicted in Figure \ref{dfsex}. The resulting order is \mbox{$b < a < c < e < d$}.

\fig{[width=\textwidth]{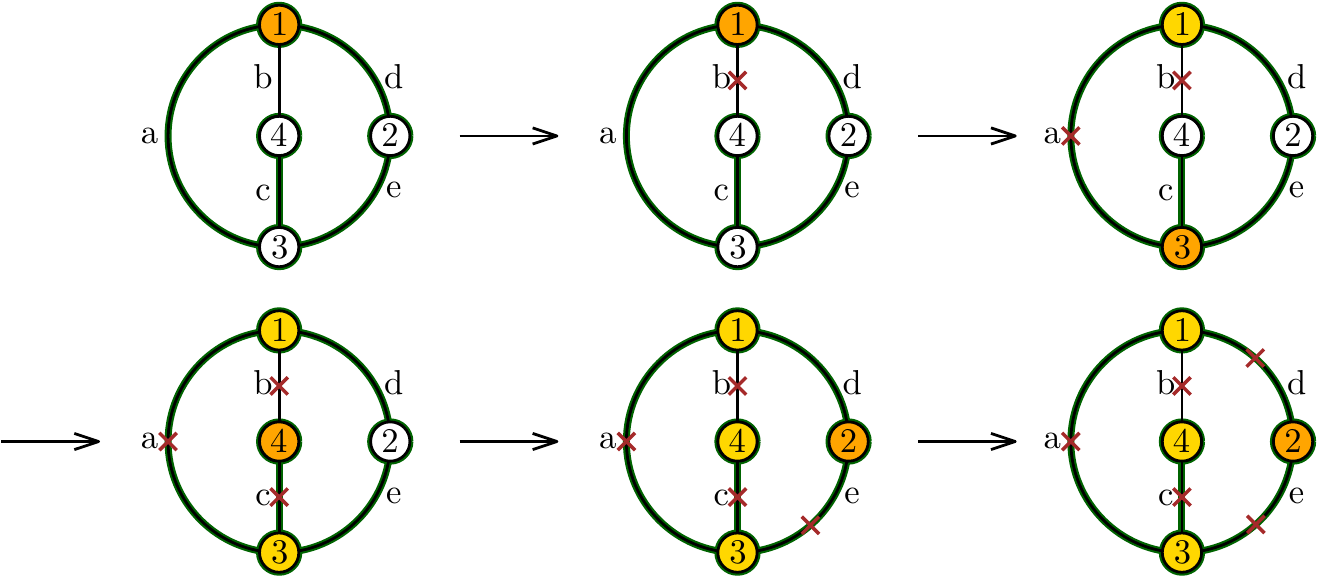}}{Illustration of a run of Algorithm \ref{pmDFS}.}{dfsex}

Let us denote by $\phi(H)$ the result of Algorithm \ref{pmDFS} for a  subgraph $H$. One can straightforwardly see that the restriction of $\phi$ to the spanning trees satisfies the hypotheses of Corollary \ref{cor:act}: indeed, the only influence of the input lies in Line \ref{lineuse}, where we test if the successive values of $\phi_j(T)$ belong to $T$ or not. 

Let $T$ be a spanning tree and $e$ an external edge. Let us prove that $e$ is DFS-active if and only if $e$ is maximal for $\phi(T)$ in its fundamental cycle $C$. If we manage to do so, we use Corollary \ref{cor:act} and the proof is ended.

a. Assume that $e$ is maximal in $C$. Since the only difference between $T$ and $T \cup e$ is $e$ (!), the executions in Algorithm \ref{pmDFS} with input $T$ and $T \cup e$ are strictly identical until the visit of $e$. In particular, at this moment, as every edge in $C$ other than $e$ has been visited and belongs the DFS forest of $T$ (because internal), both endpoints of $e$ are visited. 
Since Algorithm \ref{pmDFS} is an enriched version of Algorithm \ref{DFS}, the endpoints of $e$ are also visited just before the first visit of $e$ for Algorithm \ref{DFS} with input $T \cup e$. We never add to a DFS forest an edge whose both endpoints are marked, hence $e \notin \mathcal F(T \cup e)$.
 But $\mathcal F(T \cup e)$ is a spanning tree included in $T \cup e$. So we must have $\mathcal F(T \cup e)=T$, which means that $e$ is DFS-active.



b. Observe that for any spanning tree $\hat T$, an external edge $\hat e$ is maximal for $\phi(\hat T)$ in its fundamental cycle in $\hat T$ if and only if $\hat e$ is maximal for $\phi(\hat T \cup \hat e)$ in its fundamental cycle in $\hat T$. Indeed, the executions in Algorithm \ref{pmDFS} with input $\hat T$ and $\hat T \cup \hat e$ are the same until the visit of $\hat e$. In particular, the set of edges visited before $\hat e$ are identical.

c. Assume that $e$ is not maximal in $C$ for $\phi(T)$. 
By point b, the edge $e$ is not maximal in $C$ for $\phi(T \cup e)$.
Let $e'$ be the maximal edge in $C$ for $\phi(T \cup e)$ and denote by $T'$ the spanning tree $(T \cup e) \backslash e'$. By point b, $e'$ is maximal in $C$ for $\phi(T')$. Then, by point a, it means that $e'$ is DFS-active for $T'$. Hence $$\mathcal F(T \cup e) = \mathcal F(T' \cup e') = T' \neq T,$$
the edge $e$ is not DFS-active. 
\end{proof}

Let us conclude this section by some remarks. We have just proved that the restriction of the external DFS activity to spanning trees can be extended into a Tutte-descriptive activity. But in the light of the equations \eqref{DFStut} and \eqref{eqfor}, a question is looming: does there exist a decision tree $\Delta$ such that for all spanning forest (not necessarily spanning trees) the external DFS-active edges coincide with the edges with $\Delta$-type \bL? The answer is generally negative: consider for instance the graph  \begin{tabular}{c}\includegraphics[scale=0.8]{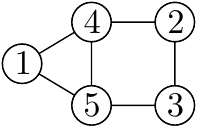}\end{tabular}. The only external DFS-active edge for the spanning forest $\ens{\ens{2,4},\ens{4,5},\ens{3,5}}$ is $\ens{2,3}$. But one can check by inspection that this edge is never DFS-active for any spanning tree. This cannot occur with $\Delta$-activities since every subgraph shares the same partition  into types as some spanning tree.

However, with some slight modifications of the definition of DFS activity, we can ensure that the correspondance between DFS-active edges and edges with type \bL\, is effective. For this, without proving it, we just have to change the input of Algorithm \ref{DFS} -- taking subgraphs of $G$ instead of general graphs -- and choose a more convenient vertex $v$ at Line \ref{linedfs}\footnote{Let us describe in a few words how to modify Line \ref{linedfs} from Algorithm \ref{DFS}. If it is the first iteration, meaning that no vertex has been visited yet, we choose $v$ as the least vertex. If some vertices have been visited, we consider $\Delta_\phi$ the decision tree induced by Corollary \ref{cor:act} and order map $\phi$, output of Algorithm \ref{pmDFS} (see proof of Proposition \ref{dfstut}). Then we run Algorithm \ref{type} with the same subgraph and with decision tree $\Delta_\phi$. We consider the edge of type \bI\, with a visited endpoint and an unvisited endpoint that was visited first. The new vertex $v$ is the unvisited endpoint of this edge.}. This new definition does not change the results of Gessel and Sagan about the DFS activity \cite{GesselSagan}, like Proposition \ref{lemGS}.

Furthermore, the external DFS activity can derive from another notion of $\Delta$-activity, named $\Delta$-forest activity. This will be the subject of Section \ref{ss:partial}. 


\section{Final comments}
\label{s:com}
In this section, we make a few comments about $\Delta$-activity and mention some prospects.

\subsection{Generalization to matroids}

The generalization of the $\Delta$-activity to the matroids is rather immediate.
Indeed, as the notions of contraction, deletion, cycle, cocycle, isthmus, loop also exist in the world of matroids, the $\Delta$-activity can be defined for matroids \textit{verbatim} (except we no longer speak about ``edges" but more generally about ``elements"). Moreover, it is not difficult to check that all the results from Sections  \ref{s:alg}, \ref{s:ord} and \ref{s:partition} hold in the same manner.
Let us add to this a property of duality, the proof of which is rather straighforward:

\begin{prop}
Let $M$ and $M^*$ be two dual matroids and $\Delta$ a decision tree. We define the decision tree $\Delta^\#$ as the mirror of $\Delta$\footnote{In term of decision functions, it means that $\Delta^\#(d_1,\dots,d_k) = \Delta(\overline d_1,\dots,\overline d_k)$, where $\overline \ell = r$ and $\overline r=\ell$.}. Given any subset $S$ of $M$, an edge is $\Delta$-active for $S$ in $M$ if and only if it is $\Delta^\#$-active for $S$ in $M^*$.
\end{prop}

\noindent \textbf{Remark. }We have not introduced the $\Delta$-activities directly on matroids because the notions of activities that we wanted to unify are more based on graphs than matroids. Moreover, I think that everyone who is familiar with matroids can generalize the notion of $\Delta$-activity without any difficulty. But the converse is not particularly true: those who are used to study the Tutte polynomial on graphs (or maps) could have been discouraged if this part was written in terms of matroids.

\subsection{Forest activities}
\label{ss:partial}
The example of DFS activity suggests that a broader notion of external activity can be adapted for spanning forests. This is the topic of this subsection. For the purpose of brevity, no proof will be given, but they can be easily adapted from the rest of this paper. 

Let $G$ be a graph. A \textit{forest activity} is a function that maps a spanning forest $F$ onto a subset of $\mathcal G(F)$, where $\mathcal G(F)$ is the set of external edges $e$ such that $F \cup e$ has a cycle. Every such edge has a \textit{fundamental cycle}, which is the unique cycle included in $F \cup e$. A forest activity $\epsilon$ is \textit{Tutte-descriptive} if the Tutte polynomial of $G$ equals
\begin{equation}
T_G(x,y) = \sum_{F\textrm{ spanning forest of }G} x^{\cc(F)-1} y^{|\epsilon(F)|}.
\end{equation}
The external DFS activity is an example of Tutte-descriptive forest activity (see Equation ~\ref{DFStut}).

Consider now a decision tree $\Delta$. Given any subgraph $S$, Algorithm \ref{ftype} outputs a subset $\epsilon(S)$ of edges, called the set of \textit{$\Delta$-forest active} edges of $S$.

\begin{algorithm}[h!]
\caption{Computing the set of $\Delta$-forest active edges.}
\label{ftype}
\begin{algorithmic}[5]
\Require $S$ subgraph of $G$.
\Ensure A subset $\epsilon(S)$ of $E(G)$.
\State $m \leftarrow $ number of edges in $G$; 
$\epsilon(S) \leftarrow \emptyset$; 
\State $n \leftarrow$ root of $\Delta$; 
\State $H \leftarrow G$;
\For {$k$ from $1$ to $m$}
	\State $e_{k} \leftarrow$ label of $n$;
	\If{$e_k$ is not a loop in $H$ \textbf{and} $e_k \notin S$ }
			\State $n \leftarrow$ left child of $n$;
	\EndIf
	\If{$e_k$ is not a loop in $H$ \textbf{and} $e_k \in S$}
		\State {$H \leftarrow \contract H e$;} 			    
		\State $n \leftarrow$ right child of $n$;
	\EndIf 
	\If{$e_k$ is a loop in $H$}
			\State \textbf{add} $e_k$ in $\epsilon(S)$;
			\State $n \leftarrow$ left child of $n$;
	\EndIf 
\EndFor
\State \Return $\epsilon(S)$
\end{algorithmic}
\end{algorithm}

\noindent \textbf{Informal description.} We start from the edge that labels the root of $\Delta$. If this edge is external or a loop, we go to the left subtree of $\Delta$. If this edge is internal and not a loop, the edge is contracted and we go to the right subtree of $\Delta$. We repeat the process until the graph has no more edge.  An external edge is $\Delta$-forest active if it is a loop  when it is  visited. Figure \ref{schema} illustrates this description.

\fig{[scale=1.4]{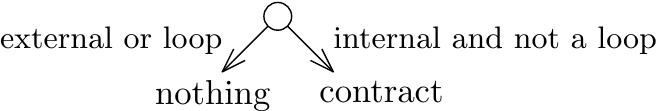}}{Diagram representing a step of Algorithm \ref{ftype}.}{schema2}

We can prove that the function that maps a spanning forest $F$ onto its set of $\Delta$-forest active edges is a Tutte-descriptive forest activity. We call it the \textit{$\Delta$-forest activity}. For any decision tree $\Delta$ the following properties are also true.
\begin{itemize}
	\item Let $F$ be a spanning forest. An edge $e$ is $\Delta$-forest active if and only if $e$ is an external edge such that $F \cup e$ has a cycle \emph{and} $e$ is maximal in its fundamental cycle for the ordering induced by the edges $e_k$ in Algorithm \ref{ftype}.
	\item The poset of subgraphs can be partitioned into subgraph intervals indexed by the spanning forests:
$$ 2^{E(G)} = \biguplus_{F\textrm{ spanning forest of }G} [F, F \cup \epsilon(F)], $$
where $\epsilon$ denotes the $\Delta$-forest activity.
  \item There exists a decision tree $\Delta'$ such that the external DFS activity coincides with the $\Delta'$-forest activity.
\end{itemize}

We could also define a dual notion of internal activity, with similar properties, based on connected subgraphs.

\subsection{Strongly Tutte-descriptive activities}

It is illusory to want to characterize every Tutte-descriptive activity. Most of them does not have any nice structure -- for example, an isthmus or a loop could be non active. We thereby need to add a constraint in order to narrow the set of interesting activities.

An activity $\psi$ is said to be \emph{strongly Tutte-descriptive} if it is Tutte-descriptive\footnote{This condition is dispensable. Indeed, any activity that satisfies \eqref{crapoeq} is Tutte-descriptive. The author thanks Vincent Delacroix for this remark.} and it induces a partition of the subgraphs:
\begin{equation} 2^{E(G)} = \biguplus_{T\textrm{ spanning tree of }G} [T \backslash \psi(T), T \cup \psi(T)]
\label{crapoeq}
 \end{equation}
(this equation is equivalent to \eqref{eqpart}).
We conjecture that the strongly Tutte-descriptive activities are precisely the $\Delta$-activities.

\begin{conjecture}
For any graph $G$, an activity $\psi$ is strongly Tutte-descriptive if and only if there exists a decision tree $\Delta$ such that $\psi$ equals the $\Delta$-activity.
\end{conjecture}

The right-to-left implication has been already proven in this paper through Theorem~\ref{charact} and Theorem \ref{partition}.  Furthermore, the above conjecture is equivalent to the following one.
\begin{conjecture}
Let $G$ be a graph with a standard edge and $\psi$ a strongly Tutte-descriptive activity. There exists an edge $e$ of $G$ such that for every spanning tree $T$, the edge $e$ does not belong to $\psi(T)$ (i.e. $e$ is active in no spanning tree).
\end{conjecture}

I will be very glad if a proof or a counter-example is found for this conjecture. Here is a proof of the equivalence between the two conjectures.

\begin{proof}
\textbf{Conjecture 1 implies Conjecture 2.}  Let us prove by induction on the total number of loops and isthmuses  of $G$ that for each decision tree $\Delta$, there exists an edge $e$ that is  $\Delta$-active in no spanning tree. Let $e$ be the label of the root node of $\Delta$. If $e$ is standard (which happens when there is no isthmus nor loop in $G$), then $e$ cannot have type \bL\, nor \bI. Consequently, $e$ cannot be $\Delta$-active for any spanning tree. If $e$ is a loop (resp. an isthmus), then $e$ is external (resp. internal) in any spanning tree. So it will be deleted (resp. contracted) in the first step in Algorithm \ref{type} and we go to the left part (resp. right part) of $\Delta$, denoted by $\Delta'$. We then use the induction hypothesis with graph $\delete G e$ (resp. $\contract G e$) and decision tree $\Delta'$. 

 \textbf{Conjecture 2 implies Conjecture 1.} Let $\psi$ be a strongly Tutte-descriptive activity for a graph $G$. 
  First we have to prove that loops and isthmuses are active for any spanning tree $T$. Let $e$ be a loop and set $S = T \cup e$. By definition of a strongly Tutte-descriptive activity, there exists a spanning tree $T'$ such that $S \in \, [T' \backslash \psi(T'),T' \cup \psi(T')]$. Since $e$ is a loop, we have $e \notin T'$. But $e \in S$, hence $e \in \psi(T')$. The spanning tree $T = S \backslash e$ must thereby belong to $[T' \backslash \psi(T'),T' \cup \psi(T')]$. But the union in \eqref{crapoeq} is disjoint, so $T=T'$. Thus $e \in \psi(T)$. We can similarly show that every isthmus belongs to $\psi(T)$.
 
    Now let us build by induction (on the number of edges) a decision tree $\Delta$ such that $\psi$ equals the $\Delta$-activity. If $G$ has standard edges, we label the root of $\Delta$ by an edge $e$ that is active in no spanning tree. (We have used Conjecture 2.) The edge $e$ is standard since it is not active for $\psi$. So $e$ cannot be $\Delta$-active since it labels the root node of $\Delta$.
Then we build a strongly Tutte-descriptive activity $\psi_c$ for $\contract G e$ by setting $\psi_c(T) = \psi(T \cup e)$ and a strongly Tutte-descriptive activity $\psi_d$ for $\delete G e$ by setting $\psi_d(T) = \psi(T)$. (The proof that these activities are indeed strongly Tutte-descriptive are left to the reader.) By induction, there exist two decision trees $\Delta_c$ et $\Delta_d$ that correspond to $\psi_c$ and $\psi_d$. The left subtree of $\Delta$ is taken to be $\Delta_d$ and the right one  is taken to be $\Delta_c$. This $\Delta$-activity is equal to $\psi$ since every spanning tree $T$ with $e$ external satisfies $\psi(T) = \psi_d(T)$ and every spanning tree $T$ with $e$ internal satisfies $\phi(T) = \psi_c(T \backslash e)$.

There remains the case (which includes the base case of the induction) where $G$ has only loops and isthmuses. However, as any loop and  any isthmus is active for $\psi$, any decision tree fits.         
\end{proof}

\section*{Acknowledgements}

The author is very grateful to his adviser Mireille Bousquet-M\'elou for her unfailing support and helpful feedbacks on this work. He would also like to thank Emeric Gioan for his pertinent advices, Robert Cori for the interim replacement of Mireille and Vincent Delacroix for his very interesting remarks. Final thanks go to Olivier Bernardi and his thesis manuscript \cite{bernathese} that was very useful for the writing of this paper.

\bibliographystyle{alpha} 
\bibliography{coloured}

\begin{thebibliography}{EMM11}

\bibitem[Ber06]{bernathese}
O.~Bernardi.
\newblock {\em Combinatoire des cartes et polyn\^ome de {T}utte}.
\newblock PhD thesis, 2006.

\bibitem[Ber08]{bernardi-tutte}
O.~Bernardi.
\newblock A characterization of the {T}utte polynomial via combinatorial
  embedding.
\newblock {\em Annals of {C}ombinatorics}, 12(2):139--153, 2008.

\bibitem[BFG02]{bdg2002}
J.~Bouttier, P.~Di Francesco, and E.~Guitter.
\newblock Census of planar maps: from the one-matrix model solution to a
  combinatorial proof.
\newblock {\em Nuclear Physics B}, 645(3):477 -- 499, 2002.

\bibitem[BO92]{broxley}
T.~Brylawski and J.~Oxley.
\newblock The {T}utte polynomial and its applications.
\newblock In {\em Matroid applications}, volume~40 of {\em Encyclopedia Math.
  Appl.}, pages 123--225. Cambridge Univ. Press, Cambridge, 1992.

\bibitem[CM92]{cori-machi}
{R}. {C}ori and {A}. {M}achi.
\newblock Maps, hypermaps and their automorphisms : a survey. {I}, {I}{I},
  {I}{I}{I}.
\newblock {\em {E}xposition. {M}ath.}, 10(5):403--467, 1992.

\bibitem[Cor75]{cori-these}
R.~Cori.
\newblock {\em Un code pour les graphes planaires et ses applications}.
\newblock Soci\'et\'e Math\'ematique de France, Paris, 1975.
\newblock With an English abstract, Ast{\'e}risque, No. 27.

\bibitem[Cou14]{courtiel-these}
J.~Courtiel.
\newblock {\em {Combinatorics of the Tutte polynomial and planar maps}}.
\newblock Theses, {Universit{\'e} de Bordeaux}, October 2014.

\bibitem[Cra69]{crapo}
H.~H. Crapo.
\newblock The {T}utte polynomial.
\newblock {\em Aequationes {M}ath.}, 3:211--229, 1969.

\bibitem[EMM11]{chap-tutte}
J.~A. Ellis-Monaghan and C.~Merino.
\newblock Graph polynomials and their applications {I}: {T}he {T}utte
  polynomial.
\newblock In {\em Structural analysis of complex networks}, pages 219--255.
  Birkh\"auser/Springer, New York, 2011.

\bibitem[FK72]{fk}
C.~M. Fortuin and P.~W. Kasteleyn.
\newblock On the random cluster model: {I}. {I}ntroduction and relation to
  other models.
\newblock {\em Physica}, 57:536--564, 1972.

\bibitem[Gre76]{curtis}
C.~Greene.
\newblock Weight enumeration and the geometry of linear codes.
\newblock {\em Studies in {A}ppl. {M}ath.}, 55(2):119--128, 1976.

\bibitem[GS96]{GesselSagan}
{I}. Gessel and {B}. Sagan.
\newblock The {T}utte polynomial of a graph, depth-first search, and simplicial
  complex partitions.
\newblock {\em The Electronic Journal of Combinatorics}, 3(2), 1996.

\bibitem[GT90]{gordon-traldi}
G.~Gordon and L.~Traldi.
\newblock Generalized activities and the {T}utte polynomial.
\newblock {\em Discrete {M}ath.}, 85(2):167--176, 1990.

\bibitem[OW79]{oxley-welsh}
J.~G. Oxley and D.~J.~A. Welsh.
\newblock The {T}utte polynomial and percolation.
\newblock In {\em Graph theory and related topics ({P}roc. {C}onf., {U}niv.
  {W}aterloo, {W}aterloo, {O}nt., 1977)}, pages 329--339. Academic Press, New
  York, 1979.

\bibitem[Sch97]{Sch97}
G.~Schaeffer.
\newblock Bijective census and random generation of {E}ulerian planar maps with
  prescribed vertex degrees.
\newblock {\em Electron. J. Combin.}, 4(1):Research Paper 20, 14 pp.\
  (electronic), 1997.

\bibitem[Thi87]{jonespoly}
M.~B. Thistlethwaite.
\newblock A spanning tree expansion of the {J}ones polynomial.
\newblock {\em Topology}, 26(3):297--309, 1987.

\bibitem[Tut54]{tutte54}
W.~T. Tutte.
\newblock A contribution on the theory of chromatic polynomial.
\newblock {\em {C}anadian {J}ournal of {M}athematics}, 6:80--91, 1954.

\bibitem[Whi32]{whitney}
H.~Whitney.
\newblock {A logical expansion in mathematics}.
\newblock {\em Bull. Am. Math. Soc.}, 38:572--579, 1932.

\end{thebibliography}

\end{document}